\documentclass[10pt]{article}

\usepackage{amssymb}
\usepackage{amsmath}
\usepackage{amsbsy}
\usepackage{amsthm}

\theoremstyle{plain}
\newtheorem{theo}{Theorem}[section]
\newtheorem{lem}{Lemma}[section]
\newtheorem{co}{Corollary}[section]
\newtheorem{pro}{Proposition}[section]
\newtheorem{rmk}{Remark}[section]
\newtheorem{defn}{Definition}[section]

\newcommand{\supp}[1]{{\rm supp\/}(#1)}
\newcommand{\sign}[1]{{\rm sign\/}(#1)}
\newcommand{\sg}[1]{{\rm sg\/}(#1)}
\begin{document}

\begin{center}
{\bf RELATIVE ASYMPTOTIC OF MULTIPLE ORTHOGONAL POLYNOMIALS FOR
NIKISHIN SYSTEMS}\\
{\bf Abey L\'opez Garc\'ia} \footnotemark \footnotetext{Dpto. de
Matem\'{a}ticas, Universidad Carlos III de Madrid, Avda.
Universidad 15, 28911 Legan\'{e}s, Madrid (SPAIN)
$<$ablopez@math.uc3m.es$>$.} \hspace{1cm} {\bf Guillermo L\'{o}pez
Lagomasino} \footnotemark \footnotetext{Dpto. de Matem\'{a}ticas,
Universidad Carlos III de Madrid, Avda. Universidad 15, 28911
Legan\'{e}s, Madrid (SPAIN) $<$lago@math.uc3m.es$>$.}
\end{center}

\begin{abstract}
We prove relative asymptotic for the ratio of two  sequences of
multiple orthogonal polynomials with respect to Nikishin system of
measures. The first Nikishin system
${\mathcal{N}}(\sigma_1,\ldots,\sigma_m)$ is such that for each
$k$, $\sigma_k$ has constant sign on its compact support   $\supp
{\sigma_k} \subset \mathbb{R}$ consisting of an interval
$\widetilde{\Delta}_k$, on which $|\sigma_k^{\prime}|
> 0$ almost everywhere, and a  discrete set without accumulation points in
$\mathbb{R} \setminus \widetilde{\Delta}_k$. If $\mbox{Co}(\supp
{\sigma_k}) = \Delta_k$ denotes the smallest interval containing
$\supp {\sigma_k}$, we assume that $\Delta_k \cap \Delta_{k+1} =
\emptyset$, $k=1,\ldots,m-1$. The second Nikishin system
${\mathcal{N}}(r_1\sigma_1,\ldots,r_m\sigma_m)$ is a perturbation
of the first by means of rational functions $r_k$, $k=1,\ldots,m,$
whose zeros and poles  lie in $\mathbb{C} \setminus \cup_{k=1}^m
\Delta_k$.
\end{abstract}



\section{Introduction}
Let $\sigma_1, \sigma_2$ be two finite Borel measures, whose
supports $\supp{\sigma_1},$ $ \supp{\sigma_2}$ are contained in
non intersecting intervals $\Delta_1,\Delta_2,$ respectively, of
the real line ${\mathbb{R}}$. Set
\[ d \langle \sigma_1, \sigma_2 \rangle(x) = \int \frac{d
\sigma_2(t)}{x-t}\, d \sigma_1(x) = \widehat{\sigma}_2(x)d
\sigma_1(x) \,.
\]
This expression defines a new measure  whose support coincides
with that of $\sigma_1$. Whenever we find it convenient we use the
differential notation of a measure.

Let $\Sigma = (\sigma_1,\ldots,\sigma_m)$ be a system of finite
Borel measures on the real line with  compact support.  $\Delta_k$
denotes the smallest interval  containing the support of
$\sigma_k$. Assume that $\Delta_k\cap\Delta_{k+1}=\emptyset,\;
k=1,\ldots,m-1.$ By definition (see \cite{kn:Nikishin}), $S =
(s_1,\ldots,s_m)= {\mathcal{N}}(\Sigma)$ is called the {\em
Nikishin system} generated by $\Sigma$ if
\[
s_1 = \sigma_1, \quad s_2 = \langle \sigma_1,\sigma_2 \rangle,
\ldots , \quad s_m = \langle \sigma_1, \langle
\sigma_2,\ldots,\sigma_m \rangle \rangle = \langle \sigma_1,
\sigma_2,\ldots,\sigma_m \rangle \,.
\]

In the sequel, $(\sigma_1,\ldots,\sigma_m)$ will always denote a
system of measures such that for each $k=1,\ldots,m, \sigma_k$ has
constant sign on its support (the sign may depend on $k$). We will
write
$(s_1,\ldots,s_m)={\mathcal{N}}^*(\sigma_1,\ldots,\sigma_m)$, if
additionally for each $k=1,\ldots,m,$  $\supp {\sigma_k} \subset
\mathbb{R}$ consists of an interval $\widetilde{\Delta}_k$, on
which $|\sigma_k^{\prime}|
> 0$ almost everywhere, and a  discrete set without accumulation points in
$\mathbb{R} \setminus \widetilde{\Delta}_k$. Finally,
$(\widetilde{s}_1,\ldots,\widetilde{s}_m)={\mathcal{N}}(p_1\sigma_1,\ldots,p_m\sigma_m),$
denotes a Nikishin system where the $p_k, k=1,\ldots,m,$ are monic
polynomials with complex coefficients whose zeros lie in
$\mathbb{C} \setminus \cup_{k=1}^m \Delta_k$.

Let $(s_1,\ldots,s_m)= {\mathcal{N}}(\sigma_1,\ldots,\sigma_m)$
and $Q_{\mathbf{n}}$ (resp. $\widetilde{Q}_{\mathbf{n}}$) be the
monic polynomial of smallest degree (not identically equal to
zero) such that
\begin{equation} \label{eq:10}
0 = \int x^{\nu} Q_{\mathbf{n}}(x) ds_k(x)\,, \quad
\nu=0,\ldots,n_k -1\,, \quad k=1,\ldots,m\,,
\end{equation}
\begin{equation}\label{eq:ortqntilde}
0 = \int x^{\nu} \widetilde{Q}_{\mathbf{n}}(x)
d\widetilde{s}_k(x)\,, \quad \nu=0,\ldots,n_k -1\,, \quad
k=1,\ldots,m\,,
\end{equation}
where $\mathbf{n} = (n_1,\ldots,n_m) \in \mathbb{Z}_+^m$. Set
$|\mathbf{n}|=n_1+\cdots+n_{m}$.

In \cite{AG} (see also \cite{AptLopRoc}), we studied the ratio
asymptotic of sequences of multiple orthogonal polynomials with
respect to a Nikishin system of measures with constant sign
extending to this setting the Rakhmanov-Denisov theorem on ratio
asymptotic of orthogonal polynomials on the real line (see
\cite{Den}, \cite{kn:Rak3}, and Proposition  \ref{cor3} below).
Here, we find the asymptotic behavior of sequences formed by
quotients of the form $\widetilde{Q}_{\bf n}/{Q_{\bf n}}$.

Given the collection of polynomials $(p_1,\ldots,p_m)$, we define
\[  \mathbb{Z}_+^m(\circledast;p_1,\ldots,p_m) = \{{\bf n}   \in \mathbb{Z}_+^m:  j <
k \Rightarrow n_k + \deg(p_{j+1}\cdots p_k) \leq n_j +1\}\,.
\]
In particular,
\[ \mathbb{Z}_+^m(\circledast) = \{{\bf n}   \in \mathbb{Z}_+^m:  j <
k \Rightarrow n_k  \leq n_j +1\}\,.
\]

A point $z_0 \in \mathbb{C}$ is said to be a $1$ attraction point
of zeros of a sequence of functions $\{\varphi_{\bf n}\}, {\bf n}
\in \Lambda \subset \mathbb{Z}_+^m,$ if for each sufficiently
small $\varepsilon
> 0$ there exists $N$ such that for all ${\bf n} \in \Lambda, |{\bf
n}| >N,$ the number of zeros (counting multiplicity) of
$\varphi_{\bf n}$  in $\{z: |z-z_0| < \varepsilon\}$ is $1$. A set
$E$ is an attractor of the zeros of $\{\varphi_{\bf n}\}, {\bf n}
\in  \Lambda,$ if for each $\varepsilon > 0$ there exists $N_0$
such that $|\mathbf{n}|>N_0, \mathbf{n}\in \Lambda,$ implies that
all the zeros of $\varphi_{{\bf n}}$ lie in the $\varepsilon$
neighborhood of $E$. Our main result states:

\begin{theo} \label{teo:1}
Let $S = {\mathcal{N}}^*(\sigma_1,\ldots,\sigma_m)$ and $\Lambda
\subset {\mathbb{Z}}_+^{m}(\circledast;p_1,\ldots,p_m)$ be a
sequence of distinct multi-indices such that for all ${\bf n} \in
\Lambda, n_{1} - n_{m} \leq  C,$ where $C$ is a constant. Then
\begin{equation} \label{eq:lim}
\lim_{{\bf n  \in \Lambda}} \frac{\widetilde{Q}_{{\bf
n}}(z)}{{Q}_{{\bf n}}(z)} = {\mathcal{F}}(z;p_1,\ldots,p_m)\,,
\quad K \subset \overline{\mathbb{C}} \setminus \supp{\sigma_1}
\,,
\end{equation}
uniformly on each compact subset $K$ of $\overline{\mathbb{C}}
\setminus \supp {\sigma_1} $, where $ {\mathcal{F}}$ is analytic
and never vanishes in $\overline{\mathbb{C}} \setminus
\widetilde{\Delta}_1$. For all sufficiently large $|\mathbf{n}|,
\mathbf{n} \in \Lambda$, $\deg \widetilde{Q}_{{\bf n}} =
|\mathbf{n}|$, $\supp{\sigma_1}$ is an attractor of the zeros of
$\{\widetilde{Q}_{\bf n}\}, {\bf n} \in \Lambda$, and each point
in $\supp{\sigma_1} \setminus \widetilde{\Delta}_1$ is a $1$
attraction point of zeros of $\{\widetilde{Q}_{\bf n}\}, {\bf n}
\in \Lambda$. When the coefficients of the polynomials $p_k,
k=1,\ldots,m,$ are real,  the statements remain valid for $\Lambda
\subset \mathbb{Z}_+^m(\circledast).$
\end{theo}

An expression for $ {\mathcal{F}}(z;p_1,\ldots,p_m)$ is given in
(\ref{expresion}) at the end of the proof of Theorem \ref{teo:1}
in section 4 below. In the sequel, any limit following the
notation used in (\ref{eq:lim}) stands for uniform convergence on
each compact subset of the indicated region.

The paper is organized as follows. Sections 2 and 3 contain
auxiliary results needed for the proof of Theorem \ref{teo:1}.
Section 4 is dedicated to its proof and deriving several
consequences. For example, we show that the same result is valid
if we modify the measures in the initial system by rational
functions instead of polynomials. These results allow to extend
the Rakhmanov-Denisov theorem on ratio asymptotic to the sequence
$\{\widetilde{Q}_{\bf n}\}, {\bf n} \in \Lambda$. In sections 5
and 6 we study the relative asymptotic of an associated system of
second type functions and their zeros.

\section{Some lemmas}

Obviously, $\mathbb{Z}_+^m(\circledast;p_1,\ldots,p_m) \subset
\mathbb{Z}_+^m(\circledast)$. If ${\bf n} \in
\mathbb{Z}_+^m(\circledast)$, it is well known  that there exists
a unique polynomial $Q_{\bf n}$ of degree $\leq |\mathbf{n}|$
satisfying the orthogonality relations expressed in (\ref{eq:10}).
Moreover, $Q_{\mathbf{n}}$ has exactly $|\mathbf{n}|$ simple zeros
which lie in the interior of $\Delta_1$ (for example, see
\cite{kn:Gonchar}).

Let us express the orthogonality relations (\ref{eq:ortqntilde})
satisfied by the polynomials $\widetilde{Q}_{\bf n}$ in terms of
the measures in the initial system.

\begin{lem} \label{lem:1} For each $k=1,\ldots,m$, we have
\begin{equation} \label{eq:3} \widetilde{s}_k = p_1 l_{k,1} s_1 +
p_1p_2 l_{k,2} s_2 + \cdots + (p_1\cdots p_k) l_{k,k} s_k\,,
\end{equation}
where $l_{k,j}$ is a polynomial of degree $\deg l_{k,j} \leq  \deg
(p_{j+1}\cdots p_k)  - 1, j < k,$ and $l_{k,k} \equiv 1.$ In
particular, if $\mathbf{n} \in
\mathbb{Z}_+^m(\circledast;p_1,\ldots,p_m)$, then
\begin{equation} \label{eq:4} 0 =
\int x^{\nu} \widetilde{Q}_{\mathbf{n}}(x) (p_1\cdots p_k)(x)
ds_k(x) \,, \quad \nu=0,\ldots,n_k-1\,,\quad k=1,\ldots,m\,.
\end{equation}
\end{lem}

{\bf Proof.} To prove (\ref{eq:3}), we proceed by induction on
$m$, the number of measures which generate the system. For $m=1$,
(\ref{eq:3}) is trivial, since  $\widetilde{s}_1 = p_1\sigma_1 =
p_1 s_1.$ Assume that (\ref{eq:3}) is true for any Nikishin system
with
 $m-1 \geq 1$ generating measures and let us prove it when the number
of generating measures is $m$.

Fix $k \in \{1,\ldots,m\}$. By definition,
\[ \widetilde{s}_k = \langle p_1\sigma_1,\ldots, p_{k} \sigma_{k}\rangle =
\langle p_1\sigma_1,\langle p_2\sigma_2, \ldots, p_{k}
\sigma_{k}\rangle \rangle\,.
\]
Consider the Nikishin system $  \mathcal{N}(p_2\sigma_2,\ldots,
p_k\sigma_k)$ which has at most $m-1$ generating measures. By the
induction hypothesis, there exist polynomials $h_{2},\ldots,
h_{k},$ $ \deg h_{j} \leq  \deg (p_{j+1}\cdots p_k)   -1, h_k
\equiv 1,$ such that
\[ \langle p_2 \sigma_2,\ldots,p_k \sigma_k\rangle = p_2
h_2 \sigma_2 + \cdots + (p_2\cdots p_k) h_k \langle
\sigma_2,\ldots,\sigma_k \rangle\,.
\]
Inserting this relation above, we have
\begin{equation} \label{eq:5}
\widetilde{s}_k = \langle p_1\sigma_1, p_2 h_2 \sigma_2 \rangle +
\cdots + \langle p_1\sigma_1, (p_2\cdots p_k) h_k \langle
\sigma_2,\ldots,\sigma_k \rangle \rangle\,.
\end{equation}

Given two measures $\sigma_{\alpha}, \sigma_{\beta}$, and a
polynomial $h$, notice that
\[ d \langle \sigma_{\alpha}, h \sigma_{\beta} \rangle (x)= \int
\frac{(h(t) \mp h(x)) d \sigma_{\beta}(t)}{x-t}
d\sigma_{\alpha}(x) =\] \[  {h}^*(x) d \sigma_{\alpha}(x) + h(x) d
\langle \sigma_{\alpha}, \sigma_{\beta} \rangle ( x)\,,
\]
where $\deg  {h}^* = \deg h -1\,.$ Making use of this property in
each term of (\ref{eq:5}), it follows that
\[ \widetilde{s}_k =
p_1[(p_2h_2)^*+\cdots + (p_2\cdots p_k h_k)^*] \sigma_1 +
(p_1p_2)h_2 \langle \sigma_1,\sigma_2 \rangle + \cdots
\]
\[    + (p_1\cdots
p_k)h_k \langle \sigma_1,\ldots,\sigma_k \rangle\,,
\]
which establishes (\ref{eq:3}).

Using (\ref{eq:3}) and the orthogonality relations
(\ref{eq:ortqntilde}) satisfied by $\widetilde{Q}_{\mathbf{n}}$,
it follows that for each $k \in \{1,\ldots,m\}$ and $\nu =
0,\ldots,n_k -1,$
\begin{equation}
\label{eq:6} 0 = \int x^{\nu} \widetilde{Q}_{\mathbf{n}}(x)
d\widetilde{s}_k(x) = \sum_{j=1}^k \int x^{\nu}l_{k,j}(x)
\widetilde{Q}_{\mathbf{n}}(x) (p_1\cdots p_j)(x) ds_j(x)\,.
\end{equation}

In the rest of the proof we assume that $\mathbf{n} \in
\mathbb{Z}_+^m(\circledast;p_1,\ldots,p_m)$. When $k=1$ the last
formula reduces to (\ref{eq:4}). Suppose that (\ref{eq:4}) holds
up to $k -1, 1 \leq k-1 \leq m-1,$ and let us show that it is also
satisfied for $k$.

Let $j \in \{1,\ldots,k-1\}$ and $0 \leq \nu \leq n_k -1$, then
\[ \nu + \deg l_{k,j} \leq n_k -1 + \deg(p_{j+1}\cdots p_k) -1 \leq
n_j -1\,.
\]
Therefore, according to the induction hypothesis
\[ \int x^{\nu}l_{k,j}(x)
\widetilde{Q}_{\mathbf{n}}(x) (p_1\cdots p_j)(x) ds_j(x) = 0\,,
\]
and (\ref{eq:6}) reduces to (\ref{eq:4}) for the index $k$. With
this we conclude the proof. \hfill $\Box$

\begin{lem} \label{lem:2} Let
$\mathbf{n} \in \mathbb{Z}_+^m(\circledast;p_1,\ldots,p_m)$. Then,
for each $k=1,\ldots,m$,
\begin{equation} \label{eq:7} 0 =
\int x^{\nu} \widetilde{Q}_{\mathbf{n}}(x) (p_1\cdots p_m)(x)
ds_k(x) \,, \,\,\, \nu=0,\ldots,n_k - \deg(p_{k+1}\cdots p_m)-1
\,.
\end{equation}
\end{lem}

{\bf Proof.} In place of $x^{\nu}$ we can put in (\ref{eq:4}) any
polynomial of degree $\leq n_k -1$. So, replacing $x^{\nu}$ by
$x^{\nu}(p_{k+1}\cdots p_m)$ we obtain (\ref{eq:7}). \hfill $\Box$

Our next objective is to express the  multiple orthogonal
polynomials of the perturbed system in terms of multiple
orthogonal polynomials of the initial system.

Let $\mathbf{n} \in \mathbb{Z}_+^m(\circledast;p_1,\ldots,p_m)$
and consider the multi-indices
\[ \mathbf{n}_j = (n_1 - \deg(p_2\cdots p_m) +j, n_2 - \deg (p_3\cdots p_m),\ldots,
n_m)\,, \quad j \geq 0\,.
\]
It is easy to verify that
\[ \mathbf{n}_j \in \mathbb{Z}_+^m(\circledast)\,, \qquad j \geq
0\,.
\]
Therefore, $\deg Q_{\mathbf{n}_j} = |\mathbf{n}_j|$ and all its
$|\mathbf{n}_j|$ simple zeros lie on $\Delta_1$. Moreover, for
each $j \geq 0$ and $k=1,\ldots,m,$
\begin{equation} \label{eq:9} 0=\int x^{\nu} Q_{\mathbf{n}_j}(x)
ds_k(x) \,, \,\,\, \nu=0,\ldots,n_k - \deg(p_{k+1}\cdots p_m)-1
\,.
\end{equation}

\begin{lem} \label{lem:3}
Let $\mathbf{n} \in \mathbb{Z}_+^m(\circledast;p_1,\ldots,p_m)$
and set $R_{\mathbf{n}} = \widetilde{Q}_{\mathbf{n}}p_1\cdots
p_m$. There exist unique constants $\lambda_{\mathbf{n},j} , j=
0,\ldots, N,$ such that
\begin{equation} \label{eq:8}
R_{\mathbf{n}} = \sum_{j=0}^{N} \lambda_{\mathbf{n},j}
Q_{\mathbf{n}_j}\,, \qquad N = \deg(p_1p_2^2\cdots p_m^m)\,.
\end{equation}
If $j'$ is such that $\deg R_{\mathbf{n}} = \deg
Q_{\mathbf{n}_{j'}}$ then  $\lambda_{\mathbf{n},j'} = 1$ and
$\lambda_{\mathbf{n},j}=0, j'+1 \leq j \leq N$. In particular,
$\lambda_{\mathbf{n},N} = 1$ if and only if $\deg
\widetilde{Q}_{\mathbf{n}} = |\mathbf{n}|$.
\end{lem}

{\bf Proof.} Since $\deg R_{\mathbf{n}} \leq |\mathbf{n}| +
\deg(p_1\cdots p_m)$, and $\{Q_{\mathbf{n}_j}\}, j = 0,\ldots ,
N,$ has representatives of all degrees from $|\mathbf{n}| -
\deg(p_2p_3^2\cdots p_m^{m-1})$ up to $|\mathbf{n}| + \deg
(p_1\cdots p_m)$, there exists a unique system of constants
$\lambda_{\mathbf{n},j} , j= 0\ldots,N,$ such that
\[
\deg (R_{\mathbf{n}} - \sum_{j=0}^{N} \lambda_{\mathbf{n},j}
Q_{\mathbf{n}_j}) \leq |\mathbf{n}| - \deg(p_2p_3^2\cdots
p_m^{m-1}) -1\,.
\]
From (\ref{eq:7})-(\ref{eq:9}) it follows that
\[ R_{\mathbf{n}} - \sum_{j=0}^{N} \lambda_{\mathbf{n},j}
Q_{\mathbf{n}_j} \equiv 0\,,
\]
which is (\ref{eq:8}). The rest of the statements follow because
$R_{\mathbf{n}}$ is monic. \hfill $\Box$

Let  ${\bf n} \in \mathbb{Z}_+^m(\circledast;p_1,\ldots,p_m).$
Define recursively the functions
\begin{equation}\label{defnRnk}
R_{{\bf n},0}(z)=R_{{\bf n}}(z),\quad R_{{\bf n},k}(z)=\int
\frac{R_{{\bf n},k-1}(x)}{z-x} d\sigma_k(x),\quad k=1,\ldots,m.
\end{equation}
In deriving  (\ref{eq:7}), we lost some orthogonality relations.
We will recover them  in the form of analytic properties of the
functions $R_{{\bf n},k}, k=0,\ldots,m-1.$

\begin{lem} \label{lem:4} Fix $\mathbf{n} \in
\mathbb{Z}_+^m(\circledast;p_1,\ldots,p_m).$ The following
relations take place:

\noindent If $z_1$ is a zero of $p_1\cdots p_m$ of multiplicity
$\tau_1,$ then
\begin{equation} \label{eq:11} \Omega_{\mathbf{n}}^{(i)}(z_1) =
\left(\frac{R_{\mathbf{n}}}{Q_{\mathbf{n}_0}}\right)^{(i)}(z_1) =
0\,, \qquad i= 0,\ldots,\tau_1 -1\,.
\end{equation}
If $z_k$ is a zero of $p_k\cdots p_m, k =2,\ldots,m,$ of
multiplicity $\tau_k,$ then
\begin{equation} \label{eq:12}
R_{\mathbf{n},k-1}^{(i)}(z_k)= 0 \,, \qquad i= 0,\ldots, \tau_k
-1\,.
\end{equation}
\end{lem}

{\bf Proof.} The zeros of $p_1\cdots p_m$ lie in $\mathbb{C}
\setminus \Delta_1$, and those of $Q_{\mathbf{n}_0}$  in
$\Delta_1$. Therefore, $\Omega_{\mathbf{n}}$ has a zero at $z_1$
of multiplicity greater or equal to $\tau_1$ which implies
(\ref{eq:11}).

For simplicity, first we will prove (\ref{eq:12})  for $k=2$. By
definition
\[ R_{{\bf n},1}(z)=\int \frac{R_{{\bf
n}}(x)}{z-x} d\sigma_1(x)\,.
\]
Therefore, for each $i \geq 0$,
\[ R_{{\bf n},1}^{(i)}(z)=(-1)^i i! \int \frac{R_{{\bf
n}}(x)}{(z-x)^{i+1}} d\sigma_1(x)\,, \qquad z \in \mathbb{C}
\setminus \Delta_1\,.
\]
If $z_2$ is a zero of $p_2\cdots p_m$ of multiplicity $\tau_2$,
using (\ref{eq:4}) with $k=1$ we have that
\[ 0 = \int \frac{(p_2\cdots p_m)(x)}{(z_2 - x)^{i+1}}
\widetilde{Q}_{\mathbf{n}}(x) p_1(x) d\sigma_1(x) =
\frac{(-1)^{i}R_{\mathbf{n},1}^{(i)}(z_2)}{i!}\,, \quad
i=0,\ldots, \tau_2-1\,,
\]
which is (\ref{eq:12}) for $k=2$. The proof of the general case
uses basically the same arguments.

Consider the functions
\[ \Phi_{\mathbf{n},k} (z) = \int \frac{R_{{\bf
n}}(x)}{z-x} ds_k(x) \,,\qquad k =1,\ldots,m\,.
\]
Notice that $\Phi_{\mathbf{n},1} = R_{\mathbf{n},1} $. For each
$i\geq 0$,
\[\Phi_{\mathbf{n},k}^{(i)} (z) = (-1)^i i! \int \frac{R_{{\bf
n}}(x)}{(z-x)^{i+1}} ds_k(x) \,,\qquad k =1,\ldots,m\,.
\]
It is easy to verify that for each $k=2,\ldots,m$,
\[\Phi_{\mathbf{n},k} (z) + (-1)^k R_{\mathbf{n},k} (z) =
\int \cdots \int \frac{R_{\mathbf{n}}(x_1)(x_1 -
x_k)d\sigma_1(x_1)\cdots
d\sigma_k(x_k)}{(z-x_1)(x_1-x_2)\cdots(x_{k-1}-x_k)(z-x_k)}\,.
\]
Taking $x_1 - x_k = x_1 - x_2 + \cdots + x_{k-1} -x_k$, it follows
that
\begin{equation} \label{eq:13} R_{\mathbf{n},k} (z) =
(-1)^{k-1} \Phi_{\mathbf{n},k} (z) + \sum_{l=1}^{k-1}(-1)^{l-1}
\widehat{\vartheta}_{l,k}(z) \Phi_{\mathbf{n},l} (z),\quad z \in
\mathbb{C} \setminus \left( \cup_{l=1}^m \Delta_l \right),
\end{equation}
where $ \vartheta_{l,k}   = \langle
\sigma_k,\sigma_{k-1},\ldots,\sigma_{l+1} \rangle$. If $z_k$ is a
zero of $p_k\cdots p_m$ of multiplicity $\tau_k (\leq \tau_{k-1}
\leq \cdots  \leq \tau_2)$, using (\ref{eq:4}) we obtain that for
each $l = 2,\ldots,k$ and $i=0,\ldots, \tau_k-1$,
\begin{equation} \label{eq:14} 0 =
\int \frac{(p_l\cdots p_m)(x)}{(z_k - x)^{i+1}}
\widetilde{Q}_{\mathbf{n}}(x) (p_1\cdots p_{l-1})(x) ds_{l-1}(x) =
\frac{(-1)^{i} {\Phi}_{\mathbf{n},l-1}^{(i)}(z_k)}{i!} \,.
\end{equation}
Now, (\ref{eq:12}) is a consequence of (\ref{eq:13}) (with $k$
replaced by $k-1$), and (\ref{eq:14}). With this we conclude the
proof. \hfill $\Box$

\section{Known asymptotic properties}
 For each ${\bf n} \in \mathbb{Z}_+^m(\circledast),$ define
recursively the   functions
\begin{equation}\label{defnPsink}
\Psi_{{\bf n},0}(z)=Q_{{\bf n}}(z),\quad\Psi_{{\bf n},k}(z)=\int
\frac{\Psi_{{\bf n},k-1}(x)}{z-x} d\sigma_k(x),\quad k=1,\ldots,m.
\end{equation}
In Proposition 1 of \cite{kn:Gonchar} it was proved that for each
${\bf n}=(n_1 ,\ldots,n_m)\in \mathbb{Z}_+^m(\circledast),$
$k=1,\ldots,m,$ and $k\leq k+r\leq m$,
\[
\int \Psi_{{\bf n},k-1}(t)\,t^\nu
d\langle\sigma_k,\ldots,\sigma_{k+r}\rangle (t)=0,\quad
\nu=0,\ldots,n_{k+r}-1.
\]
From here, the authors deduce that $\Psi_{{\bf n},k-1}, k =
1,\ldots,m,$ has exactly $N_{{\bf n},k} = n_{k}+\cdots+n_m$ zeros
in ${\mathbb{C}} \setminus \Delta_{k-1}$, that they are all
simple, and lie in the interior of $\Delta_{k}$. Let $Q_{{\bf
n},k}$ be the monic polynomial of degree $N_{{\bf n},k}$ whose
simple zeros are located at the points where $\Psi_{{\bf n},k-1}$
vanishes on $\Delta_{k}$ and let $Q_{{\bf n},m+1} \equiv 1$. In
Proposition 2 (see also Proposition 3) of \cite{kn:Gonchar} the
authors show that
\begin{equation}               \label{eq:defineqn}
\int x^{\nu} \Psi_{{\bf n},k-1}(x)\frac{d\sigma_k(x)}{Q_{{\bf
n},k+1}(x)}=0,\quad \nu=0,\ldots,N_{{\bf n},k}-1,\quad
k=1,\ldots,m\,.
\end{equation}
Set
\[H_{{\bf n},k}(z) :=  \frac{Q_{{\bf n},k-1}(z)\Psi_{{\bf n},k-1}(z)}{Q_{{\bf
n},k}(z)}\,, \qquad k=1,\ldots,m+1\,,
\]
$(H_{{\bf n},1}(z)\equiv 1).$ It is well known (see (50) in
\cite{kn:B-G}) and easy to verify that
\begin{equation} \label{eq:18} H_{{\bf
n},k+1}(z) = \int \frac{Q_{\mathbf{n},k}^2(x)}{z-x} \frac{H_{{\bf
n},k}(x)d\sigma_{k}(x)}{Q_{\mathbf{n},k-1}(x)Q_{\mathbf{n},k+1}(x)}\,,\qquad
k = 1,\ldots,m\,.
\end{equation}

From (\ref{eq:defineqn}), we have that for each multi-index ${\bf
n}=(n_1,\ldots,n_m)\in \mathbb{Z}_+^m(\circledast)$ there exists
an associated system of polynomials
\[
\{Q_{{\bf n},k}\}_{k=1}^m ,\quad \deg Q_{{\bf
n},k}=\sum_{\alpha=k}^m n_\alpha=:N_{{\bf n},k} ,\qquad Q_{{\bf
n},0}\equiv Q_{{\bf n},m+1}\equiv 1.
\]
For each $k=1,\ldots,m,$ they satisfy the full system of
orthogonality relation
\begin{equation}            \label{eq:ortogonalvar}
\int x^{\nu} Q_{{\bf n},k}(x)\frac{H_{{\bf
n},k}(x)d\sigma_k(x)}{Q_{{\bf n},k-1}(x) Q_{{\bf
n},k+1}(x)}=0,\quad \nu=0,\ldots,N_{{\bf n},k}-1,
\end{equation}
with respect to varying measures. Notice that  $H_{{\bf n},k}$ and
$Q_{{\bf n},k-1} Q_{{\bf n},k+1}$ have constant sign on
$\Delta_k$.

Let $\varepsilon_{\mathbf{n},k}$ be the sign of the measure $
{H_{{\bf n},k}(x)d\sigma_k(x)}/{Q_{{\bf n},k-1}(x) Q_{{\bf
n},k+1}(x)}$ on $\supp{\sigma_k}$. For each $k=1,\ldots,m,$ set
\begin{equation} \label{eq:K} K_{{\bf n},k} = \left( \int
Q_{{\bf n},k}^2(x)
\frac{\varepsilon_{\mathbf{n},k}H_{\mathbf{n},k}(x)d\sigma_k(x)}{Q_{{\bf
n},k-1}(x)Q_{{\bf n},k+1}(x)} \right)^{-1/2}\;.
\end{equation}
Take
\[ K_{{\bf n},0} = 1 \;, \quad \kappa_{{\bf n},k} = \frac{K_{{\bf n},k}}{K_{{\bf n},k-1}}
\;, \quad k=1,\ldots,m \;.
\]
Define
\begin{equation} \label{eq:orton} q_{{\bf n},k} =
\kappa_{{\bf n},k}Q_{{\bf n},k} \;, \quad h_{{\bf n},k}  = K_{{\bf
n},k-1}^2 H_{{\bf n},k} \;, \quad k = 1,\ldots,m \;.
\end{equation}
From (\ref{eq:ortogonalvar})
\[
\int x^{\nu} Q_{{\bf
n},k}(x)\frac{\varepsilon_{\mathbf{n},k}h_{{\bf
n},k}(x)d\sigma_k(x)}{Q_{{\bf n},k-1}(x) Q_{{\bf
n},k+1}(x)}=0,\quad \nu=0,\ldots,N_{{\bf n},k}-1,\quad
k=1,\ldots,m\,,
\]
and with the notation introduced above it follows that $q_{n,k}$
is orthonormal with respect to the varying measure
\[
\frac{\varepsilon_{\mathbf{n},k}h_{{\bf
n},k}(x)d\sigma_k(x)}{Q_{{\bf n},k-1}(x)Q_{{\bf n},k+1}(x)} =
d\rho_{{\bf n},k}(x)\,.
\]

In Lemma 3.3 of \cite{AG} (see also Corollary 3 in \cite{DolBer})
we proved

\begin{pro} \label{pro1}
Let $S = {\mathcal{N}}^*(\sigma_1,\ldots,\sigma_m)$ and $\Lambda
\subset {\mathbb{Z}}_+^{m}(\circledast)$ be a sequence of
multi-indices such that for all ${\bf n} \in \Lambda,$  $n_{1} -
n_{m} \leq  C,$ where $C$ is a constant.  Then, for each fixed $k
= 1,\ldots,m,$ we have
\begin{equation} \label{eq:x}
\lim_{{\bf n}\in {\Lambda}}
\varepsilon_{\mathbf{n},k}h_{\mathbf{n},k+1}(z) =
\frac{1}{\sqrt{(z - b_{k})(z- a_{k})}}\,, \qquad K \subset
\overline{\mathbb{C}} \setminus \supp{\sigma_k}\,,
\end{equation}
where $[a_{k},b_{k}]=\widetilde{\Delta}_k$. The square root is
taken so that $\sqrt{(z-b_{k})(z-a_{k})}>0$ for $z=x>b_{k}$.
$\supp{\sigma_k}$ is an attractor of the zeros of $\{Q_{{\bf
n},k}\}, {\bf n} \in \Lambda,$ and each point of $\supp{\sigma_k}
\setminus \widetilde{\Delta}_k$ is a $1$ attraction point of zeros
of $\{Q_{{\bf n},k}\}, {\bf n} \in \Lambda$.
\end{pro}

In the proof of our main result, we use the asymptotic behavior of
the polynomials $Q_{\mathbf{n},k}, k=1,\ldots,m,$ and the
functions $\Psi_{\mathbf{n},k}, k=1,\ldots,m,$ when $\mathbf{n}$
runs through a sequence of multi-indices $\Lambda \subset
\mathbf{Z}_+^m(\circledast)$. In order to describe these
asymptotic formulas we need to introduce some notions.

Consider the $(m+1)$-sheeted Riemann surface
$$
\mathcal R=\overline{\bigcup_{k=0}^m \mathcal R_k} ,
$$
formed by the consecutively ``glued'' sheets
$$
\mathcal R_0:=\overline {\mathbb{C}} \setminus
\widetilde{\Delta}_1,\quad \mathcal R_k:=\overline {\mathbb{C}}
\setminus \{\widetilde{\Delta}_k \cup
\widetilde{\Delta}_{k+1}\},\qquad k=1,\ldots,m-1,\qquad \mathcal
R_m=\overline {\mathbb{C}} \setminus \widetilde{\Delta}_m,
$$
where the upper and lower banks of the slits of two neighboring
sheets are identified. Fix $l \in \{1,\ldots,m\}$. Let
$\psi^{(l)}, l=1,\ldots,m,$ be a single valued rational function
on $\mathcal{R}$ whose divisor consists of a simple zero at the
point $\infty^{(0)} \in \mathcal R_0$ and a simple pole at the
point $\infty^{(l)} \in \mathcal R_l$. Therefore,
\[
\psi^{(l)}(z) = C_1/z + \mathcal{O}(1/z^2)\,,\,\,z \to
\infty^{(0)}\,, \qquad  \psi^{(l)}(z) = C_2 z +
\mathcal{O}(1)\,,\,\,z \to \infty^{(l)} \,,
\]
where $C_1$ and $C_2$ are constants different from zero.  Since
the genus of $\mathcal R$ equals zero, such a single valued
function on $\mathcal R$ exists and it is uniquely determined
except for a multiplicative constant. We denote the branches of
the algebraic function $\psi^{(l)}$, corresponding to the
different sheets $k = 0,\ldots,m$ of $\mathcal R$ by
$$
\psi^{(l)}:=\{\psi_k^{(l)}\}_{k=0}^m\,.
$$
We normalize $\psi^{(l)}$ so that
\begin{equation} \label{eq:normaliz*}
\prod_{k=0}^{m}\,|\psi^{(l)}_{k}(\infty)|=1, \qquad C_1 \in
\mathbb{R}\setminus \{0\}.
\end{equation}
Certainly, there are two $\psi^{(l)}$ verifying this
normalization.

Since the product of all the branches $\prod_{k=0}^{m}
\psi^{(l)}_{k}$ is a single valued analytic function in
$\overline{\mathbb{C}}$ without singularities, by Liouville's
Theorem it is constant and because of the normalization introduced
above this constant is either $1$ or $-1$. Since $\psi^{(l)}$ is
such that $C_1 \in \mathbb{R}\setminus \{0\},$ then
\[ \psi^{(l)}(z) = \overline{\psi^{(l)}(\overline{z})}, \qquad z \in
\mathcal{R}.
\]
In fact, let $\phi(z) := \overline{\psi^{(l)}(\overline{z})}$.
$\phi$ and $\psi^{(l)}$ have the same divisor; consequently, there
exists a constant $C$ such that $\phi= C\psi^{(l)}$. Comparing the
leading coefficients of the Laurent expansion of these functions
at $\infty^{(0)}$,  we conclude that $C=1$.

Given an arbitrary function $F(z)$  which has in a neighborhood of
infinity a Laurent expansion of the form $F(z)= Cz^k +
{\mathcal{O}}(z^{k-1}), C \neq 0,$ and $ k \in {\mathbb{Z}},$ we
denote
\[
\widetilde{F}:= {F}/{C}\,.
\]
(For simplicity in writing, we write $\widetilde{F}_k^{(l)}$
instead of the more appropriate $\widetilde{F_k^{(l)}}$.) $C$ is
called the leading coefficient of $F$. When $C \in \mathbb{R}
\setminus \{0\}$, $\mbox{sg}(F(\infty))$ represents the sign of
$C$.

In terms of the branches of $\psi^{(l)}$, the symmetry formula
above means that that for each  $k= 0,1,\ldots,m$
\begin{equation} \label{real}
\psi^{(l)}_k: \overline{\mathbb{R}} \setminus
(\widetilde{\Delta}_k\cup \widetilde{\Delta}_{k+1})
\longrightarrow \overline{\mathbb{R}}
\end{equation}
$(\widetilde{\Delta}_0 = \widetilde{\Delta}_{m+1}=\emptyset)$;
therefore, the coefficients (in particular, the leading one) of
the Laurent expansion at $\infty$ of these branches are real
numbers and $\mbox{sg}(\psi^{(l)}_k(\infty))$ is  defined. It also
expresses that
\[
\psi^{(l)}_k(x_{\pm}) = \overline{\psi^{(l)}_k(x_{\mp})} =
\overline{\psi^{(l)}_{k+1}(x_{\pm})}, \qquad x \in
\widetilde{\Delta}_{k+1}.
\]

For any fixed multi-index ${\bf n}=(n_1,\ldots,n_m)$, set
$$
{\bf n}^{l}:=(n_1,\ldots,n_{l-1},n_l+1,n_{l+1},\ldots,n_m)\,.
$$
In \cite{AG} (see also \cite{AptLopRoc}) the authors prove

\begin{pro} \label{cor3}
Let $S = {\mathcal{N}}^*(\sigma_1,\ldots,\sigma_m)$ and $\Lambda
\subset \mathbb{Z}^{m}_+(\circledast)$ be a sequence of
multi-indices such that for all ${\bf n} \in \Lambda$ and some
fixed $l \in \{1,\ldots,m\}$, we have that ${\bf n}^{l} \in
\mathbb{Z}^{m}_+(\circledast)$ and $n_{1} - n_{m} \leq  C,$ where
$C$ is a constant. Let $\{q_{{\bf n},k} = \kappa_{{\bf
n},k}Q_{{\bf n},k} \}_{k=1}^m ,{\bf n}\in {\Lambda},$ be the
system of orthonormal polynomials defined in $(\ref{eq:orton})$
and $\{K_{{\bf n},k}\}_{k=1}^m ,{\bf n}\in {\Lambda},$ the values
given by $(\ref{eq:K})$. Then, for each fixed $k = 1,\ldots,m,$ we
have
\begin{equation} \label{eq:xa}\lim_{{\bf n}\in
{\Lambda}}\frac{\kappa_{{\bf n}^{l},k}}{\kappa_{{\bf n},k}}=
\kappa^{(l)}_k\,,
\end{equation}
\begin{equation} \label{eq:xk}
\lim_{{\bf n}\in {\Lambda}}\frac{K_{{\bf n}^{l},k}}{K_{{\bf
n},k}}= \kappa^{(l)}_1\cdots\kappa^{(l)}_k\,,
\end{equation}
and
\begin{equation} \label{eq:xb}
\lim_{{\bf n}\in {\Lambda}}\frac{q_{{\bf n}^{l},k}(z)}{q_{{\bf
n},k}(z)}= \kappa^{(l)}_k \widetilde{F}_k^{(l)}(z), \qquad K
\subset \mathbb{C} \setminus \supp{\sigma_k} \,,
\end{equation}
where
\begin{equation} \label{eq:xc} \kappa^{(l)}_k =
\frac{c_{k}^{(l)}}{\sqrt{c_{k-1}^{(l)}c_{k+1}^{(l)}}}\,, \qquad
c_{k}^{(l)} = \left\{
\begin{array}{ll}
(F^{(l)}_k)^{\prime}(\infty)\,, & k=1,\ldots,l \,, \\
F^{(l)}_k(\infty)\,, & k= l+1,\ldots,m\,,
\end{array}
\right.
\end{equation}
$(c_{0}^{(l)}=c_{m+1}^{(l)}=1)$ and
 \begin{equation} \label{eq:F}
F_k^{(l)} := \delta_{k,l} {\prod_{\nu=k}^m \psi_{\nu}^{(l)}},
\end{equation}
with $\delta_{k,l}= \mbox{sg}\left(\prod_{\nu=k}^m
\psi_{\nu}^{(l)}(\infty)\right).$
\end{pro}

\section{Proof of Theorem \ref{teo:1}}
When $l=1$, it is possible to find an algebraic function
$\psi^{(1)}$ satisfying
\begin{equation} \label{eq:normaliz}
\prod_{k=0}^m \psi_k^{(1)}(\infty) = 1\,, \qquad
C_{1}\in\mathbb{R}\setminus\{0\}.
\end{equation}
Let $(a,b)_k$ denote the interval $(a,b)$ on the sheet
$\mathcal{R}_k$. We distinguish two cases. Suppose that
$\widetilde{\Delta}_{1}=[a_1,b_1]$ is to the left of
$\widetilde{\Delta}_{2}=[a_2,b_2]$. Take $\psi^{(1)}$ verifying
(\ref{eq:normaliz*}) with
$C_1=\lim_{z\rightarrow\infty}z\psi_{0}^{(1)}(z)>0$. Because of
(\ref{real}),   the restriction of  $\psi^{(1)}$ to
$(-\infty,a_1]_0 \cup (-\infty,a_1]_1$ establishes a bicontinuous
bijection onto the interval $(-\infty,0)$ of the real line. It
follows that $\psi^{(1)}_1(x) \to -\infty, x \to -\infty, x \in
\mathbb{R},$ which means that $C_2 >0$, and $\psi^{(1)}_k(\infty)
> 0, k=2,\ldots,m.$ Therefore,
$\prod_{k=0}^m \psi_k^{(1)}(\infty) > 0.$ If
$\widetilde{\Delta}_{1}$ is to the right of
$\widetilde{\Delta}_{2}$, take $\psi^{(1)}$ satisfying
(\ref{eq:normaliz*}) with  $C_1<0$. Now, the restriction of
$\psi^{(1)}$ to $[b_1,+\infty)_0 \cup [b_1,+\infty)_1$ establishes
a bicontinuous bijection onto  $(-\infty,0)$. It follows that
$\psi^{(1)}_1(x) \to -\infty, x \to +\infty, x \in \mathbb{R},$
which means that $C_2 < 0$, and $\psi^{(1)}_k(\infty)
> 0, k=2,\ldots,m.$ Again, $\prod_{k=0}^m \psi_k^{(1)}(\infty) > 0.$

Throughout the rest of the paper, when $\widetilde{\Delta}_{1}$ is
to the left of $\widetilde{\Delta}_{2},$ we will select
$\psi^{(1)}$ so that $\mbox{sg}(\psi_{k}^{(1)}(\infty))=1$, for
all $k=0,\ldots,m$. If $\widetilde{\Delta}_{1}$ is to the right of
$\widetilde{\Delta}_{2}$, we will take $\psi^{(1)}$ so that
$\mbox{sg}(\psi_{0}^{(1)}(\infty))=\mbox{sg}(\psi_{1}^{(1)}(\infty))=-1$
and $\mbox{sg}(\psi_{k}^{(1)}(\infty))=1$, for all $k=2,\ldots,m$.

In general, for any $l \in\{1,\ldots,m\}$ and $\psi^{(l)}$
verifying (\ref{eq:normaliz*}), we know that
\[\prod_{\nu=0}^{m}\psi_{\nu}^{(l)}(\infty) \in \{-1,1\}\,.\]

Let $\Lambda \subset \mathbb{Z}_+^m(\circledast;p_1,\ldots,p_m)$
be an infinite sequence of distinct multi-indices such that $n_1 -
n_m \leq C,\mathbf{n} \in \Lambda$. According to
(\ref{eq:xa})-(\ref{eq:F}), for each fixed $j \geq 0,$
\begin{equation} \label{eq:30}\lim_{\mathbf{n} \in \Lambda} \frac{
{Q}_{\mathbf{n}_{j+1}}(z)}{ {Q}_{\mathbf{n}_j}(z)} =
\widetilde{F}_1^{(1)}(z) =
\frac{\mbox{sg}(\psi_0^{(1)}(\infty))}{c_1^{(1)}\psi_0^{(1)}(z)}
=: \varphi_0(z), \quad  K \subset {\mathbb{C}} \setminus
\supp{\sigma_1} \,.
\end{equation}
(Notice that (\ref{eq:normaliz}) implies that $\prod_{\nu= 0}^m
\psi_{\nu}^{(1)}(z) \equiv 1.)$

Using (\ref{eq:8}),
\[\Omega_{\mathbf{n}} =  \frac{R_{\mathbf{n}}}{Q_{\mathbf{n}_0}} = \sum_{j=0}^{N} \lambda_{\mathbf{n},j}
\frac{Q_{\mathbf{n}_j}}{Q_{\mathbf{n}_0}}\,, \qquad N =
\deg(p_1p_2^2\cdots p_m^m)\,.
\]
Set
\[ \lambda_{\mathbf{n}}^* = (\sum_{j=0}^N
|\lambda_{\mathbf{n},j}|)^{-1}\,.
\]
At least one of the numbers in the sum is $1$ so
$\lambda_{\mathbf{n}}^*$ is finite. Define
\begin{equation}
\label{eq:nor} \lambda_{\mathbf{n}}^*{\Omega}_{\mathbf{n}} =
\sum_{j=0}^{N} \lambda_{\mathbf{n},j}^*
\frac{Q_{\mathbf{n}_j}}{Q_{\mathbf{n}_0}}\,, \qquad \sum_{j=0}^{N}
|\lambda_{\mathbf{n},j}^*| =1\,.
\end{equation}

Because of (\ref{eq:30}) and (\ref{eq:nor}), the family $\{
\lambda_{\mathbf{n}}^*{\Omega}_{\mathbf{n}}\}, \mathbf{n} \in
\Lambda,$ is normal in $\mathbb{C} \setminus \supp{\sigma_1}$, and
any convergent subsequence $\{
\lambda_{\mathbf{n}}^*{\Omega}_{\mathbf{n}}\}, \mathbf{n} \in
\Lambda' \subset \Lambda,$ converges to
\[ \lim_{ \mathbf{n} \in \Lambda'} \lambda_{\mathbf{n}}^*{\Omega}_{\mathbf{n}}(z) =
p_{\Lambda'} ( {\varphi_0(z)} ) = \sum_{j=0}^{N}
\lambda_{j}\varphi_0^j(z)\,, \quad K \subset \mathbb{C} \setminus
\supp{\sigma_1}\,.
\]
That is, $p_{\Lambda'}(w)$ is a polynomial of degree $\leq N$, not
identically equal to zero since $\sum_{j=0}^{N} |\lambda_{j}| =
1$. We will show that $p_{\Lambda'}$ does not depend on the
subsequence taken. This implies the existence of limit along all
$\Lambda$. To this aim, we will uniquely determine $N$ zeros of
$p_{\Lambda'}$.

Let $z_1$ be  one of the zeros of $p_1\cdots p_m$ and $\tau_1$ its
multiplicity. Using (\ref{eq:11}) and the Weierstrass theorem, it
follows that
\[ (p_{\Lambda'} \circ \varphi_0)^{(i)}(z_1) = 0\,, \qquad i=0,\ldots,\tau_1
-1\,.
\]
Since $\varphi_0$ is one to one in $\mathbb{C} \setminus
\widetilde{\Delta}_1$, we conclude that $p_{\Lambda'}(w)$ is
divisible by
\[ (w - \varphi_0(z_1))^{\tau_1}\,.
\]

We will detect the rest of the zeros of $p_{\Lambda'}(w)$ in
virtue of (\ref{eq:12}). Consider the sequence
$\{\lambda_{\mathbf{n}}^* R_{\mathbf{n},k-1}\}, \mathbf{n} \in
\Lambda'$. From (\ref{eq:8}), (\ref{defnRnk}) and
(\ref{defnPsink})
\[ \lambda_{\mathbf{n}}^* R_{\mathbf{n},k-1}(z) =
\sum_{j=0}^N \lambda_{\mathbf{n},j}^*
{\Psi}_{\mathbf{n}_j,k-1}(z)\,.
\]
Multiplying this equation by
$\varepsilon_{\mathbf{n}_0,k-1}K_{\mathbf{n}_0,k-1}^2
Q_{\mathbf{n}_0,k-1}/Q_{\mathbf{n}_0,k}$ and using the definition
of $h_{\mathbf{n},k}$,  we obtain
\[
\frac{\lambda_{\mathbf{n}}^*
\varepsilon_{\mathbf{n}_0,k-1}K_{\mathbf{n}_0,k-1}^2
(Q_{\mathbf{n}_0,k-1}R_{\mathbf{n},k-1})(z)}{Q_{\mathbf{n}_0,k}(z)}
\]
\[ =
 \sum_{j=0}^N \lambda_{\mathbf{n},j}^* \frac{K_{\mathbf{n}_0,k-1}^2}{K_{\mathbf{n}_j,k-1}^2}
\frac{Q_{\mathbf{n}_0,k-1}(z)}{Q_{\mathbf{n}_j,k-1}(z)}
\frac{Q_{\mathbf{n}_j,k}(z)}{Q_{\mathbf{n}_0,k}(z)}
\frac{\varepsilon_{\mathbf{n}_0,k-1}}{\varepsilon_{\mathbf{n}_j,k-1}}
\varepsilon_{\mathbf{n}_j,k-1}h_{\mathbf{n}_j,k}(z)\,.
\]

From (\ref{eq:xa})-(\ref{eq:xb}), for each $j \geq 0$ and
$k=2,\ldots,m,$
\[ \lim_{\mathbf{n} \in \Lambda'}\frac{K_{\mathbf{n}_j,k-1}^2}{K_{\mathbf{n}_{j+1},k-1}^2}
\frac{Q_{\mathbf{n}_j,k-1}(z)}{Q_{\mathbf{n}_{j+1},k-1}(z)}
\frac{Q_{\mathbf{n}_{j+1},k}(z)}{Q_{\mathbf{n}_j,k}(z)} =
\frac{\widetilde{F}_k^{(1)}(z)}{(\kappa^{(1)}_1\cdots\kappa^{(1)}_{k-1})^2\widetilde{F}_{k-1}^{(1)}(z)}\,,
\]
uniformly on compact subsets of $\mathbb{C} \setminus
(\supp{\sigma_{k-1}} \cup \supp{\sigma_k})$. On account of
(\ref{eq:xc}) and the expression of the functions $F_k^{(1)}$,
\begin{equation}\label{eqdefvarphi}
\frac{\widetilde{F}_k^{(1)}(z)}{(\kappa^{(1)}_1\cdots\kappa^{(1)}_{k-1})^2\widetilde{F}_{k-1}^{(1)}(z)}
= \frac{\sg{
\psi_{k-1}^{(1)}(\infty)}}{c_1^{(1)}\psi_{k-1}^{(1)}(z)}
=:\varphi_{k-1}(z).
\end{equation}

Let us consider the ratios
$\varepsilon_{\mathbf{n}_{j+1},k}/\varepsilon_{\mathbf{n}_{j},k},
k= 1,\ldots, m-1, j \geq 0.$ Recall that
$\varepsilon_{\mathbf{n},k}$ is by definition the sign of the
measure
$H_{\mathbf{n},k}(x)d\sigma_k(x)/(Q_{\mathbf{n},k-1}Q_{\mathbf{n},k+1})(x)$
on $\Delta_k$. Notice that for each fixed $k=2,\ldots,m$ the
polynomials $Q_{\mathbf{n}_j,k}$ have the same degree for all $j
\geq 0$; therefore,   they all have the same sign on any interval
disjoint from $\Delta_k$. On the other hand, the polynomials
$Q_{\mathbf{n}_j,1}$ have degrees that increase one by one with
$j$. Hence, if $\Delta_1$ is to the left of $\Delta_2$, all the
polynomials $Q_{\mathbf{n}_j,1}$ have the same sign on $\Delta_2$
whereas, if $\Delta_1$ is to the right of $\Delta_2$, the sign of
these polynomials alternate on $\Delta_2$ as $j$ increases one by
one. Taking these facts into consideration, it is easy to see that
for all $j\geq 0$, the measures
$H_{\mathbf{n}_j,1}(x)d\sigma_1(x)/ Q_{\mathbf{n}_j,2}(x) =
d\sigma_1(x)/ Q_{\mathbf{n}_j,2}(x),$ have the same sign;
therefore, for all $j \geq 0$,
$\varepsilon_{\mathbf{n}_{j+1},1}/\varepsilon_{\mathbf{n}_{j},1} =
1 $ and the functions $H_{\mathbf{n}_j,2} $ have the same sign on
$\Delta_2$ (see (\ref{eq:18})). Hence, the measures
$H_{\mathbf{n}_j,2}(x)d\sigma_2(x)/
(Q_{\mathbf{n}_j,1}Q_{\mathbf{n}_j,3})(x)  $ have the same sign if
$\Delta_1$ is to the left of $\Delta_2$ and alternate signs as $j$
increases when $\Delta_1$ is to the right of $\Delta_2$. Thus, for
all $j \geq 0$,
$\varepsilon_{\mathbf{n}_{j+1},2}/\varepsilon_{\mathbf{n}_{j},2} =
1 $ when $\Delta_1$ is to the left of $\Delta_2$ and
$\varepsilon_{\mathbf{n}_{j+1},2}/\varepsilon_{\mathbf{n}_{j},2} =
-1 $ when $\Delta_1$ is to the right of $\Delta_2$. By the same
token (see (\ref{eq:18})), for all $j \geq 0$ the functions
$H_{\mathbf{n}_j,3} $ have the same sign on $\Delta_3$ when
$\Delta_1$ is to the left of $\Delta_2$ and alternate sign when
$\Delta_1$ is to the right of $\Delta_2$. From now on the
situation repeats and for each fixed $k=2,\ldots,m-1,$  and all $j
\geq 0$,
$\varepsilon_{\mathbf{n}_{j+1},k}/\varepsilon_{\mathbf{n}_{j},k} =
1 $ when $\Delta_1$ is to the left of $\Delta_2$ while
$\varepsilon_{\mathbf{n}_{j+1},k}/\varepsilon_{\mathbf{n}_{j},k} =
-1 $ when $\Delta_1$ is to the right of $\Delta_2$.

Let $\delta = 1$ when $\Delta_1$ is to the left of $\Delta_2$ and
$\delta = -1$ if $\Delta_1$ is to the right of $\Delta_2$. Using
(\ref{eq:x}) and (\ref{eq:xa})-(\ref{eq:xc}), it follows that
\[ \lim_{\mathbf{n} \in \Lambda'} \lambda_{\mathbf{n}}^*
\varepsilon_{\mathbf{n}_0,k-1}K_{\mathbf{n}_0,k-1}^2  \frac{
Q_{\mathbf{n}_0,k-1}(z)R_{\mathbf{n},k-1}(z)}{Q_{\mathbf{n}_0,k}(z)}
= \]
\begin{equation} \label{eq:32}
\left\{
\begin{array}{ll}
\frac{1}{\sqrt{(z-b_{1})(z-a_{1})}}\sum_{j=0}^N \lambda_{j}
\varphi_{1}^j(z)\,, & k=2\,, \\
\frac{1}{\sqrt{(z-b_{k-1})(z-a_{k-1})}}\sum_{j=0}^N \lambda_{j}
(\delta\varphi_{k-1})^j(z)\,, & k=3,\ldots,m\,, \\
\end{array} \right. =
\end{equation}
\[
\left\{
\begin{array}{ll}
\frac{1}{\sqrt{(z-b_{1})(z-a_{1})}} p_{\Lambda'}(\varphi_{1}(z))\,,& k=2\,, \\
\frac{1}{\sqrt{(z-b_{k-1})(z-a_{k-1})}}p_{\Lambda'}(\delta\varphi_{k-1}(z))\,, & k=3,\ldots,m\,, \\
\end{array} \right.
\]
uniformly on each compact subset $K$ of $\mathbb{C} \setminus
(\supp{\sigma_{k-1}} \cup \supp{\sigma_k})$.

Let $z_k$ be one of the zeros of $p_k\cdots p_m, k=2,\ldots,m,$
and $\tau_k$ its multiplicity. Using (\ref{eq:32}), (\ref{eq:12}),
and the Weierstrass theorem, it follows that
\[ (p_{\Lambda'} \circ \varphi_{1})^{(i)}(z_2) = 0\,, \qquad
i=0,\ldots,\tau_2 -1\,,
\]
and
\[ (p_{\Lambda'} \circ (\delta\varphi_{k-1}))^{(i)}(z_k) = 0\,, \quad
i=0,\ldots,\tau_k -1\,, \quad k =3,\ldots,m \,.
\]
Since $\varphi_{k-1}$ is one to one in $\mathbb{C} \setminus
(\widetilde{\Delta}_{k-1} \cup \widetilde{\Delta}_k)$, we conclude
that $p_{\Lambda'}(w)$ is divisible by
\[ (w - \varphi_{1}(z_2))^{\tau_2}\,,
\]
and
\[ (w - \delta\varphi_{k-1}(z_k))^{\tau_k}\,, \quad k=3,\ldots,m\,.
\]
Therefore, the following sets are formed by zeros of
$p_{\Lambda'}$:
\[\mathcal{Z}_{0}:=\{\varphi_{0}(z_1): z_1 \,\mbox{is a zero of}\, p_1\cdots p_{m}\}\,,\]
\[\mathcal{Z}_{1}:=\{\varphi_{1}(z_2): z_2 \,\mbox{is a zero of}\, p_2\cdots p_{m}\}\,,\]
\[\mathcal{Z}_{k}:=\{\delta\varphi_{k}(z_{k+1}): z_{k+1}\,\mbox{is a zero of}\,
 p_{k+1}\cdots p_{m}\}\,,\quad 2\leq k\leq m-1\,.\]

Assume first that $\delta=1$. Recall that in this case we selected
$\psi^{(1)}$ so that $\mbox{sg}(\psi^{(1)}_{k}(\infty))=1$ for all
$0\leq k\leq m$. Therefore the functions $\varphi_{0},
\varphi_{1}, \delta\varphi_{k}, 2\leq k\leq m-1,$ are the first
$m$ branches of $1/c_{1}^{(1)}\psi^{(1)}$. If $\delta=-1$, since
$\psi^{(1)}$ was chosen so that
$\mbox{sg}(\psi^{(1)}_{0}(\infty))=\mbox{sg}(\psi^{(1)}_{1}(\infty))=-1$
and $\mbox{sg}(\psi^{(1)}_{k}(\infty))=1, 2\leq k\leq m$, the
functions $\varphi_{0}, \varphi_{1}, \delta\varphi_{k}, 2\leq
k\leq m-1,$ are now the first $m$ branches of
$-1/c_{1}^{(1)}\psi^{(1)}$. In any case, since $\psi^{(1)}:
\mathcal{R}\longrightarrow \overline{\mathbb{C}}$ is bijective it
follows that the zero sets $\mathcal{Z}_{k}, 0\leq k\leq m-1$ are
pairwise disjoint. Therefore, we have detected
$N=\deg(p_1p_2^2\cdots p_m^m)$ zeros (counting multiplicities) of
the polynomial $p_{\Lambda'}$ and their location does not depend
on the subsequence $\Lambda' \subset \Lambda$.

Let
\[(p_k\cdots p_m)(z) = \prod_{\nu = 1}^{l_k} (z -
z_{k,\nu})^{\tau_{k,\nu}}\,,
\]
where $\{z_{k,1},\ldots,z_{k,l_k}\}$ are the distinct zeros of
$p_k\cdots p_m$. Then
\[ p_{\Lambda'}(w) = c \prod_{k=1}^2 \prod_{\nu = 1}^{l_k} (w -
\varphi_{k-1}(z_{k,\nu}))^{\tau_{k,\nu}}\prod_{k=3}^m \prod_{\nu =
1}^{l_k} (w - \delta\varphi_{k-1}(z_{k,\nu}))^{\tau_{k,\nu}}\,,
\]
where $c$ is uniquely defined by the conditions that it is a
positive constant such that the sum of the moduli of the
coefficients of $p_{\Lambda'}$ must equal one; moreover,
\[ 0 < c = \lim_{\mathrm{n} \in \Lambda} \lambda_{\mathbf{n}}^* <
\infty\,.
\]
Consequently, uniformly on each compact subset $K \subset
\mathbb{C} \setminus \supp{\sigma_1}$,
\[
\lim_{\mathbf{n} \in \Lambda}
\frac{R_{\mathbf{n}}(z)}{Q_{\mathbf{n}_0}(z)} =
\]
\begin{equation} \label{eq:lim1} \prod_{k=1}^2
\prod_{\nu = 1}^{l_k} (\varphi_0(z) -
\varphi_{k-1}(z_{k,\nu}))^{\tau_{k,\nu}}\prod_{k=3}^m \prod_{\nu =
1}^{l_k} (\varphi_0(z) - \delta
\varphi_{k-1}(z_{k,\nu}))^{\tau_{k,\nu}}\,.
\end{equation}

From (\ref{eq:xa}) and  (\ref{eq:xb}), it follows that
\begin{equation} \label{eq:lim3}
\lim_{\mathbf{n} \in \Lambda}
\frac{Q_{\mathbf{n}}(z)}{Q_{\mathbf{n}_0}(z)}=(\widetilde{F}_1^{(1)}(z))^{\deg(p_2
\cdots p_m)} \cdots(\widetilde{F}_1^{(m-1)}(z))^{\deg(p_m)}\,.
\end{equation}
Combining (\ref{eq:lim1}) and (\ref{eq:lim3}), we get
\[
\lim_{\mathbf{n} \in \Lambda}
\frac{\widetilde{Q}_{\mathbf{n}}(z)}{Q_{\mathbf{n}}(z)} =
 {\mathcal{F}}(z;p_1,\ldots,p_m)\,, \qquad K \subset
\mathbb{C} \setminus \supp{\sigma_1}\,,
\]
where ($\varphi_{0}(z)=\widetilde{F}_1^{(1)}(z)$)
\[{\mathcal{F}}(z;p_1,\ldots,p_m) = \prod_{\nu =1}^{l_1} \left(\frac{\varphi_0(z) -
\varphi_0(z_{1,\nu})}{z-z_{1,\nu}}\right)^{\tau_{1,\nu}}
\prod_{\nu = 1}^{l_2} \left(1 -
\frac{\varphi_1(z_{2,\nu})}{\varphi_0(z)}\right)^{\tau_{2,\nu}}
\times \] \[ \prod_{k=3}^m\prod_{\nu = 1}^{l_k} \left(
\frac{\varphi_0(z)-\delta\varphi_{k-1}(z_{k,\nu})}{\widetilde{F}_1^{(k-1)}(z)
}\right)^{\tau_{k,\nu}}\,.
\]

Let us simplify the expression above. From the definition of the
functions $\varphi_{k}$, and taking into account that
$\delta=\sg{\psi_{0}^{(1)}(\infty)}$, it follows that
\[1 - \frac{\varphi_1(z_{2,\nu})}{\varphi_0(z)}=1 -\frac{\psi_{0}^{(1)}(z)}{\psi_{1}^{(1)}(z_{2,\nu})}\,.\]

It is easy to see that for $l\geq 2$ the following equation holds:
\begin{equation}\label{eq:relalg}
\frac{1}{\psi^{(1)}(z)}-\frac{1}{\psi^{(1)}(\infty^{(l-1)})}=\frac{C_{0}^{(l-1)}}{C_{0}^{(1)}\psi^{(l-1)}(z)}\,,
\end{equation}
where
\[\psi^{(1)}(z) = C_{0}^{(1)}/z + \mathcal{O}(1/z^2)\,,\,\,z \to
\infty^{(0)}\,,\]
\[\psi^{(l-1)}(z) = C_{0}^{(l-1)}/z +
\mathcal{O}(1/z^2)\,,\,\,z \to \infty^{(0)} \,.\]

For $k\geq3$ (recall that
$\prod_{\nu=0}^{m}\psi_{\nu}^{(l)}(\infty)\in \{-1,1\}$ when
$l\geq 2$), we have that
\[\widetilde{F}_1^{(k-1)}(z)=\frac{\sg{\psi_{0}^{(k-1)}}(\infty)}{c_{1}^{(k-1)}\psi_{0}^{(k-1)}(z)}\,.\]
Thus
\begin{equation}\label{eq:neweq}
\frac{\varphi_0(z)-\delta\varphi_{k-1}(z_{k,\nu})}{\widetilde{F}_1^{(k-1)}(z)
}=\frac{c_{1}^{(k-1)}\psi_{0}^{(k-1)}(z)}{c_{1}^{(1)}\sg{\psi_{0}^{(k-1)}(\infty)}}\left(\frac{\sg{\psi_{0}^{(1)}(\infty)}}
{\psi_{0}^{(1)}(z)}-\frac{\delta}{\psi_{k-1}^{(1)}(z_{k,\nu})}\right).
\end{equation}

From (\ref{eq:relalg}), it follows that
\[\psi_{0}^{(k-1)}(z)\left(\frac{1}{\psi_{0}^{(1)}(z)}-\frac{1}{\psi_{k-1}^{(1)}(\infty)}\right)
=\frac{C_{0}^{(k-1)}}{C_{0}^{(1)}}\,.\] Therefore,
\begin{equation}\label{eq:neweq2}
\psi_{0}^{(k-1)}(z)\left(\frac{1}{\psi_{0}^{(1)}(z)}-\frac{\delta}{\psi_{k-1}^{(1)}(z_{k,\nu})}\right)=\frac{C_{0}^{(k-1)}}{C_{0}^{(1)}}
+\left(\frac{\psi_{0}^{(k-1)}(z)}{\psi_{k-1}^{(1)}(\infty)}-\frac{\delta
\psi_{0}^{(k-1)}(z)}{\psi_{k-1}^{(1)}(z_{k,\nu})}\right)\,.
\end{equation}

It is straightforward to check that
\begin{equation}\label{eq:neweq3}
\frac{c_{1}^{(k-1)}}{c_{1}^{(1)}}\frac{C_{0}^{(k-1)}}{C_{0}^{(1)}}=\frac{\sg{\psi_{0}^{(k-1)}(\infty)}}{\sg{\psi_{0}^{(1)}(\infty)}}\,.
\end{equation}
Evaluating (\ref{eq:relalg}) at $z_{k,\nu}$ we obtain
\begin{equation}\label{eq:neweq4}
\frac{1}{\psi_{k-1}^{(1)}(z_{k,\nu})}-\frac{1}{\psi_{k-1}^{(1)}(\infty)}=\frac{C_{0}^{(k-1)}}{C_{0}^{(1)}\psi_{k-1}^{(k-1)}(z_{k,\nu})}\,.
\end{equation}

Assume that $\Delta_{1}$ is to the left of $\Delta_{2}$, then
$\delta=\sg{\psi_{0}^{(1)}(\infty)}=1$. From (\ref{eq:neweq}),
(\ref{eq:neweq2}), (\ref{eq:neweq3}), and (\ref{eq:neweq4}), we
find that
\[\frac{\varphi_0(z)-\delta\varphi_{k-1}(z_{k,\nu})}{\widetilde{F}_1^{(k-1)}(z)
}=1-\frac{\psi_{0}^{(k-1)}(z)}{\psi_{k-1}^{(k-1)}(z_{k,\nu})}\,.\]
If $\Delta_{1}$ is to the right of $\Delta_{2}$, then $\delta =
\sg{\psi_{0}^{(1)}(\infty)}=-1$. Applying
(\ref{eq:neweq})-(\ref{eq:neweq4}), we obtain again
\[\frac{\varphi_0(z)-\delta\varphi_{k-1}(z_{k,\nu})}{\widetilde{F}_1^{(k-1)}(z)
}=1-\frac{\psi_{0}^{(k-1)}(z)}{\psi_{k-1}^{(k-1)}(z_{k,\nu})}\,.\]
Therefore,
\begin{equation} \label{expresion}
\mathcal{F}(z;p_1,\ldots,p_m)=\prod_{\nu =1}^{l_1}
\left(\frac{\varphi_0(z) -
\varphi_0(z_{1,\nu})}{z-z_{1,\nu}}\right)^{\tau_{1,\nu}}\prod_{k=2}^m\prod_{\nu
= 1}^{l_k}
\left(1-\frac{\psi_{0}^{(k-1)}(z)}{\psi_{k-1}^{(k-1)}(z_{k,\nu})}
\right)^{\tau_{k,\nu}}.
\end{equation}
(We did not substitute $\varphi_0$ in terms of $\psi_0^{(1)}$ (see
(\ref{eq:30})) in the first group of products for simplicity in
the final expression.)

We have proved (\ref{eq:lim}) on compact subsets of ${\mathbb{C}}
\setminus \supp{\sigma_1}$. Using the maximum principle it follows
that the same is true on compact subsets of $\overline{\mathbb{C}}
\setminus \supp{\sigma_1}$. Notice that $ {\mathcal{F}}$ is
analytic and has no zero in $\overline{\mathbb{C}} \setminus
\widetilde{\Delta}_1$. For all $\mathbf{n} \in \Lambda$, $\deg
{Q}_{{\bf n}} = |\mathbf{n}|$, $\supp{\sigma_1}$ is an attractor
of the zeros of $\{Q_{\bf n}\}, {\bf n} \in \Lambda$, and each
point in $\supp{\sigma_1} \setminus \widetilde{\Delta}_1$ is a $1$
attraction point of zeros of $\{Q_{\bf n}\}, {\bf n} \in \Lambda$;
therefore, the statements concerning $\deg
\widetilde{Q}_{\mathbf{n}}$ and the asymptotic behavior of the
zeros of these polynomials follow from (\ref{eq:lim}), on account
of the argument principle and the corresponding behavior of the
zeros of the polynomials $Q_{\bf n}$ described in Proposition
\ref{pro1}.

In order to prove the last statement, let us assume that the
polynomials $p_k, k=1,\ldots,m,$ have real coefficients and
$\Lambda \subset \mathbf{Z}_+^m(\circledast)$. Notice that in this
case the polynomials $\widetilde{Q}_{\bf n}$ are the multiple
orthogonal polynomials with respect to the Nikishin system
$\mathcal{N}(p_1\sigma_1,\ldots,p_m\sigma_m)$ generated by real
measures with constant sign. Thus, Proposition \ref{cor3} can be
applied to them. Given $\Lambda$ we construct the auxiliary
sequence $\Lambda(\diamond)$ as follows. To each ${\bf n} =
(n_1,\ldots,n_m) \in \Lambda$ we associate ${\bf n}_{\diamond} =
(n_1, n_2 - \deg(p_2),\ldots,n_m - \deg(p_2\cdots p_m))$ (we
disregard those multi-indices in $\Lambda$ for which a component
of ${\bf n}_{\diamond}$ would turn out to be negative, according
to the assumptions on $\Lambda$ there can be at most a finite
number of such $\bf n$). It is easy to see that $\Lambda(\diamond)
\subset \mathbb{Z}_+^m(\circledast;p_1,\ldots,p_m)$.

Choose consecutive multi-indices running from ${\bf n}_{\diamond}$
to $\mathbf{n}$ so that each one of them belongs to
$\mathbb{Z}_+^m(\circledast)$. We can write $Q_{\bf n}/Q_{{\bf
n}_{\diamond}}$ as the product of quotients of the corresponding
monic multiple orthogonal polynomials. The same can be done with
$\widetilde{Q}_{\bf n}/\widetilde{Q}_{{\bf n}_{\diamond}}$.
According to (\ref{eq:xa}) and (\ref{eq:xb}), there exists an
analytic function $G(z)$ in  $\mathbb{C} \setminus
\widetilde{\Delta}_1$, which is never zero, such that
\[ \lim_{{\bf n} \in \Lambda} \frac{Q_{\bf n}(z)}{Q_{{\bf
n}_{\diamond}}(z)} = \lim_{{\bf n} \in \Lambda}
\frac{\widetilde{Q}_{\bf n}(z)}{\widetilde{Q}_{{\bf
n}_{\diamond}}(z)} = G(z) \,, \qquad K \subset \mathbb{C}
\setminus \supp{\sigma_1}\,.
\]
Since
\[ \frac{\widetilde{Q}_{\bf n}(z)}{Q_{{\bf
n} }(z)} = \frac{\widetilde{Q}_{\bf n}(z)}{\widetilde{Q}_{{\bf
n}_{\diamond} }(z)} \frac{\widetilde{Q}_{{\bf
n}_{\diamond}}(z)}{Q_{{\bf n}_{\diamond}}(z)} \frac{Q_{{\bf
n}_{\diamond}}(z)}{Q_{{\bf n} }(z)}\,,
\]
using Theorem \ref{teo:1} on the ratio in the middle and the
previous limits on the other two ratios, the last statement
readily follows. \hfill $\Box$

We can easily extend the main result to the more general case when
the perturbation on the initial system is carried out by rational
functions.

\begin{co} \label{co:2}
Let $S = {\mathcal{N}}^*(\sigma_1,\ldots,\sigma_m)$. Consider the
perturbed Nikishin system $
\mathcal{N}(\frac{p_1}{q_1}\sigma_1,\ldots,\frac{p_m}{q_m}\sigma_m)$,
where $p_k,q_k$ denote relatively prime polynomials whose zeros
lie in  $\mathbb{C} \setminus \cup_{k=1}^m \Delta_k$. Let $\Lambda
\subset {\mathbb{Z}}^{m}_+(\circledast;p_1q_1,\ldots,p_mq_m)$ be a
sequence of distinct multi-indices such that for all ${\bf n} \in
\Lambda, n_{1} - n_{m} \leq  C,$ where $C$ is a constant. Let
$\widetilde{Q}_{{\bf n}}$ be the monic multiple orthogonal
polynomial of smallest degree relative to the Nikishin system
$\mathcal{N}(\frac{p_1}{q_1}\sigma_1,\ldots,\frac{p_m}{q_m}\sigma_m)$
and ${\bf n}$. Then
\begin{equation} \label{eq:lim2}
\lim_{{\bf n  \in \Lambda}} \frac{\widetilde{Q}_{{\bf
n}}(z)}{{Q}_{{\bf n}}(z)} = \frac{ {\mathcal{F}}(z; p_1 ,\ldots,
p_m)}{{\mathcal{F}}(z; q_1 ,\ldots, q_m)}\,, \quad K \subset
\overline{\mathbb{C}} \setminus \supp{\sigma_1} \,.
\end{equation}
For all sufficiently large $|\mathbf{n}|, \mathbf{n} \in \Lambda$,
$\deg \widetilde{Q}_{{\bf n}} = |\mathbf{n}|$, $\supp{\sigma_1}$
is an attractor of the zeros of $\{\widetilde{Q}_{\bf n}\}, {\bf
n} \in \Lambda$, and each point in $\supp{\sigma_1} \setminus
\widetilde{\Delta}_1$ is a $1$ attraction point of zeros of
$\{\widetilde{Q}_{\bf n}\}, {\bf n} \in \Lambda$. When the
polynomials $p_{k}, q_{k}, k=1,\ldots,m$, have real coefficients,
the statements remain valid for
$\Lambda\subset{\mathbb{Z}}^{m}_+(\circledast)$.
\end{co}

{\bf Proof.} Notice that
$\mathcal{N}(\frac{p_1}{q_1}\sigma_1,\ldots,\frac{p_m}{q_m}\sigma_m)
=
\mathcal{N}(\frac{p_1\overline{q}_1}{|q_1|^2}\sigma_1,\ldots,\frac{p_m\overline{q}_m}{|q_m|^2}\sigma_m)$,
where $\overline{q}_k$ denotes the polynomial obtained conjugating
the coefficients of $q_k$.  Let $Q_{\bf n}^*$ be the ${\bf n}$th
monic multiple orthogonal polynomial with respect to the Nikishin
system
$\mathcal{N}(\frac{\sigma_1}{|q_1|^2},\ldots,\frac{\sigma_m}{|q_m|^2})
$ generated by measures with constant sign.

Using Theorem \ref{teo:1},
\[
\lim_{{\bf n}\in \Lambda}\frac{\widetilde{Q}_{{\bf n}}(z)}{Q_{{\bf
n}}^*(z)}= {\mathcal{F}}(z;
p_1\overline{q}_1,\ldots,p_m\overline{q}_m)\,, \quad K \subset
\overline{\mathbb{C}} \setminus \supp{\sigma_1}
\]
and, considering the last remark of the same theorem, we also have
\[ \lim_{{\bf n}\in \Lambda}\frac{ {Q}_{{\bf n}}(z)}{Q_{{\bf n}}^*(z)}=
{\mathcal{F}}(z; q_1\overline{q}_1,\ldots,q_m\overline{q}_m)\,,
\quad K \subset \overline{\mathbb{C}} \setminus \supp{\sigma_1}\,.
\]
On the other hand,
\[  \frac{{\mathcal{F}}(z;
p_1\overline{q}_1,\ldots,p_m\overline{q}_m)}{{\mathcal{F}}(z;
q_1\overline{q}_1,\ldots,q_m\overline{q}_m)} =
\frac{{\mathcal{F}}(z; p_1 ,\ldots,p_m) }{{\mathcal{F}}(z; q_1
,\ldots,q_m )}
\]
because in the products defining the functions on the left hand
side all the factors connected with the zeros of the
$\overline{q}_k$ cancel out. Consequently, (\ref{eq:lim2}) takes
place.  The rest of the statements of the corollary are proved
following arguments similar to those employed in the proof of
Theorem \ref{teo:1}. \hfill $\Box$

The previous results allow to derive ratio asymptotic for the
multiple orthogonal polynomials of our perturbed Nikishin systems.

\begin{co}\label{co:3}
Let $S = {\mathcal{N}}^*(\sigma_1,\ldots,\sigma_m)$ Consider the
perturbed Nikishin system
$\mathcal{N}(\frac{p_1}{q_1}\sigma_1,\ldots,\frac{p_m}{q_m}\sigma_m)$,
where $p_k,q_k$ denote relatively prime polynomials whose zeros
lie in  $\mathbb{C} \setminus \cup_{k=1}^m \Delta_k$. Let $\Lambda
\subset {\mathbb{Z}}^{m}_+(\circledast;p_1q_1,\ldots,p_mq_m)$ be a
sequence of distinct multi-indices such that for all ${\bf n} \in
\Lambda$ and some fixed $l\in\{1,\ldots,m\}$, we have that ${\bf
n}^{l}\in {\mathbb{Z}}^{m}_+(\circledast;p_1q_1,\ldots,p_mq_m)$
and $n_{1} - n_{m} \leq C,$ where $C$ is a constant. Let
$\widetilde{Q}_{{\bf n}}$ be the monic multiple orthogonal
polynomial of smallest degree with respect to the Nikishin system
$\mathcal{N}(\frac{p_1}{q_1}\sigma_1,\ldots,\frac{p_m}{q_m}\sigma_m)$
and ${\bf n}$. Then
\[
\lim_{{\bf n}\in \Lambda}\frac{\widetilde{Q}_{{\bf
n}^{l}}(z)}{\widetilde{Q}_{{\bf n}}(z)}=\lim_{{\bf n}\in
\Lambda}\frac{Q_{{\bf n}^{l}}(z)}{Q_{{\bf
n}}(z)}=\widetilde{F}_1^{(l)}(z),\qquad K \subset \mathbb{C}
\setminus \supp{\sigma_1}\,. \]
\end{co}
{\bf Proof.} Since
\[
\frac{\widetilde{Q}_{{\bf n}^{l}}(z)}{\widetilde{Q}_{{\bf
n}}(z)}=\frac{\widetilde{Q}_{{\bf n}^{l}}(z)}{Q_{{\bf
n}^{l}}(z)}\frac{Q_{{\bf n}^{l}}(z)}{Q_{{\bf n}}(z)}\frac{Q_{{\bf
n}}(z)}{\widetilde{Q}_{{\bf n}}(z)}\,,
\]
the result follows immediately applying Proposition \ref{cor3} and
Corollary \ref{co:2}. \hfill $\Box$
\section{Relative asymptotic of second type functions}
Let $\widetilde{Q}_{{\bf n}}$ be the monic polynomial of smallest
degree satisfying (\ref{eq:ortqntilde}). Set
\[\widetilde{\Psi}_{n,0}(z):=\widetilde{Q}_{{\bf n}}(z)\,,\]
\begin{equation}\label{defnPsinktilde}
\widetilde{\Psi}_{n,k}(z):=\int\frac{\widetilde{\Psi}_{n,k-1}(x)}{z-x}p_{k}(x)\,d\sigma_{k}(x)\,,\quad
1\leq k\leq m\,.
\end{equation}
\begin{lem}\label{lemmarelation}
If $n_{j}\geq \deg(p_{j+1}\cdots p_{m})$, $j=1,\ldots,m-1$, then
$R_{{\bf n},k}(z)=(p_{k+1}\cdots p_{m})(z)\widetilde{\Psi}_{{\bf
n},k}(z)$, $z\in\mathbb{C}\setminus\supp{\sigma_{k}}$,
$k=0,1,\ldots,m,$ $(R_{{\bf n},m}=\widetilde{\Psi}_{{\bf n},m})$.
\end{lem}
{\bf Proof.} We proceed by induction on $k$. The case $k=0$ is
trivial since by definition, $R_{{\bf n},0}(z)=(p_{1}\cdots
p_{m})(z)\widetilde{Q}_{{\bf n}}(z)$. Assume that the result holds
for $k-1$, and let us prove it for $k$. We have
\[
R_{{\bf n},k}(z)=\int\frac{R_{{\bf
n},k-1}(x)}{z-x}d\sigma_{k}(x)=\int\frac{\widetilde{\Psi}_{{\bf
n},k-1}(x)(p_{k}\cdots p_{m})(x)}{z-x}d\sigma_{k}(x)=
\]
\[
(p_{k+1}\cdots p_{m})(z)\widetilde{\Psi}_{{\bf n},k}(z)+\int
\widetilde{\Psi}_{{\bf n},k-1}(x)l(x)p_{k}(x)d\sigma_{k}(x)\,,
\]
where $l(x)$ is a polynomial of degree $\deg(p_{k+1}\cdots
p_{m})-1$. Now, for $k\leq k+r \leq m,$ the functions
$\widetilde{\Psi}_{{\bf n},k}$ satisfy the orthogonality relations
(see in \cite{kn:Gonchar} that the proof exposed there is valid
also for complex measures)
\[
\int \widetilde{\Psi}_{{\bf n},k-1}(t)t^{\nu}d\langle
p_{k}\sigma_{k},\ldots,p_{k+r}\sigma_{k+r}\rangle(t)=0\,,\qquad
\nu=0,\ldots,n_{k+r}-1.
\]
In particular, $\int \widetilde{\Psi}_{{\bf
n},k-1}(t)t^{\nu}p_{k}(t)d\sigma_{k}(t)=0$ if $\nu\leq n_{k}-1$.
Thus, since we are assuming that $n_{k}\geq \deg(p_{k+1}\cdots
p_{m})$, we get that
\[
\int \widetilde{\Psi}_{{\bf n},k-1}(x)l(x)p_{k}(x)d\sigma_{k}(x)=0
\]
and the result follows. \hfill $\Box$

\begin{rmk}
The condition $n_{k}\geq \deg(p_{k+1}\cdots p_{m}),
k=1,\ldots,m-1,$ is automatically satisfied by the components of
multi-indices ${\bf n}$ with norm sufficiently large that belong
to a sequence $\Lambda \subset
{\mathbb{Z}}_+^{m}(\circledast;p_1,\ldots,p_m)$ such that for all
${\bf n} \in \Lambda, n_{1} - n_{m} \leq  C,$ where $C$ is a
constant. In fact, it is satisfied for all ${\bf n} \in
{\mathbb{Z}}_+^{m}(\circledast;p_1,\ldots,p_m)$ such that $n_m
\geq 1$.
\end{rmk}
Now, we need to introduce some notations. Let
\[
\delta_{k}:=\left\{
\begin{array}{ccccc}
1, & \mbox{if} & \Delta_{k}
& \mbox{is to the left of} & \Delta_{k+1} \,,\\
-1, & \mbox{if} & \Delta_{k}
& \mbox{is to the right of} & \Delta_{k+1} \,.\\
\end{array}
\right.
\]
For $k\geq 2$, set
\[
\Delta_{k,l}:=\left\{
\begin{array}{ccc}
-\delta_{k}\delta_{k-1}, & \mbox{if} & l\geq k+1 \,,\\
\delta_{k-1}, & \mbox{if} & l\in\{k-1, k\} \,,\\
1, & \mbox{if} & l\leq k-2 \,.
\end{array}
\right.
\]
If $k=1$,
\[
\Delta_{1,l}:=\left\{
\begin{array}{ccc}
1, & \mbox{if} & l=1 \,,\\
-\delta_{1}, & \mbox{if} & l\geq 2 \,.
\end{array}
\right.
\]
\begin{lem}
For any $\mathbf{n}, \mathbf{n}^l
\in\mathbb{Z}_{+}^{m}(\circledast)$
\begin{equation}\label{ratioepsilons}
\frac{\varepsilon_{\mathbf{n}^{l},k}}{\varepsilon_{\mathbf{n},k}}=\prod_{i=1}^{k}\Delta_{i,l}\,.
\end{equation}
\end{lem}
{\bf Proof.} By definition, $\varepsilon_{\mathbf{n},k}$ is the
sign of the measure
$\frac{H_{\mathbf{n},k}(x)d\sigma_{k}(x)}{Q_{\mathbf{n},k-1}(x)Q_{\mathbf{n},k+1}(x)}$
on $\supp{\sigma_{k}}$. We will denote by $\sign{f,\Delta}$ the
sign of a function $f$ on the interval $\Delta$. Thus
\begin{equation}\label{eqepsilon1}
\frac{\varepsilon_{\mathbf{n}^{l},k}}{\varepsilon_{\mathbf{n},k}}
=\sign{\frac{H_{\mathbf{n}^{l},k}Q_{\mathbf{n},k-1}Q_{\mathbf{n},k+1}}{
{H_{\mathbf{n},k}Q_{\mathbf{n}^{l},k-1}}
Q_{\mathbf{n}^{l},k+1}},\Delta_{k}}\,.
\end{equation}

If $l\geq k-1$, since
$\deg(Q_{\mathbf{n}^{l},k-1})=1+\deg(Q_{\mathbf{n},k-1})$, we have
that
\begin{equation}\label{eqdelta1}
\sign{ {Q_{\mathbf{n},k-1}}/{Q_{\mathbf{n}^{l},k-1}},\Delta_{k}}
=\delta_{k-1}\,,
\end{equation}
and if $l\leq k-2$, since
$\deg(Q_{\mathbf{n}^{l},k-1})=\deg(Q_{\mathbf{n}^{l},k-1})$, we
obtain
\begin{equation}\label{eqdelta2}
\sign{
{Q_{\mathbf{n},k-1}}/{Q_{\mathbf{n}^{l},k-1}},\Delta_{k}}=1\,.
\end{equation}
By similar arguments, we know that for $l\geq k+1$,
\begin{equation}\label{eqdelta3}
\sign{ {Q_{\mathbf{n},k+1}}/{Q_{\mathbf{n}^{l},k+1}},\Delta_{k}}
=-\delta_{k}\,,
\end{equation}
and if $l\leq k$,
\begin{equation}\label{eqdelta4}
\sign{
{Q_{\mathbf{n},k+1}}/{Q_{\mathbf{n}^{l},k+1}},\Delta_{k}}=1\,.
\end{equation}

Finally, from (\ref{eq:18}) it follows that
\[\frac{H_{\mathbf{n}^{l},k}(x)}{H_{\mathbf{n},k}(x)}=
\frac{\int_{\Delta_{k-1}}\frac{Q_{\mathbf{n}^{l},k-1}^{2}(t)}{x-t}\frac{H_{\mathbf{n}^{l},k-1}(t)
d\sigma_{k-1}(t)}{Q_{\mathbf{n}^{l},k-2}(t)Q_{\mathbf{n}^{l},k}(t)}}{\int_{\Delta_{k-1}}\frac{Q_{\mathbf{n},k-1}^{2}(t)}{x-t}
\frac{H_{\mathbf{n},k-1}(t)
d\sigma_{k-1}(t)}{Q_{\mathbf{n},k-2}(t)Q_{\mathbf{n},k}(t)}}\,.\]
Therefore,
\begin{equation}\label{eqepsilon2}
\sign{ {H_{\mathbf{n}^{l},k}}/{H_{\mathbf{n},k}},
\Delta_{k}}=\frac{\varepsilon_{\mathbf{n}^{l},k-1}}
{\varepsilon_{\mathbf{n},k-1}}\,.
\end{equation}

From (\ref{eqepsilon1})-(\ref{eqepsilon2}) we conclude that
\[\frac{\varepsilon_{\mathbf{n}^{l},k}}{\varepsilon_{\mathbf{n},k}}=\Delta_{k,l}\frac{\varepsilon_{\mathbf{n}^{l},k-1}}{\varepsilon_{\mathbf{n},k-1}}\,.\]
Since $H_{\mathbf{n}^{l},1} \equiv H_{\mathbf{n},1}\equiv
Q_{\mathbf{n}^{l},0} \equiv Q_{\mathbf{n},0}\equiv 1$, we have
that $\varepsilon_{\mathbf{n}^{l},1}$ is the sign of the measure
$\frac{d\sigma_{1}(x)}{Q_{\mathbf{n}^{l},2}(x)}$ on $\Delta_{1}$,
and $\varepsilon_{\mathbf{n},1}$ is the sign of the measure
$\frac{d\sigma_{1}(x)}{Q_{\mathbf{n},2}(x)}$ on $\Delta_{1}$.
Therefore, we have (\ref{ratioepsilons}). \hfill $\Box$
\begin{defn}
We define the following functions
\begin{equation}\label{defphi}
\varphi_{k-1}^{(j)}(z):=\frac{\sg{\psi_{k-1}^{(j)}(\infty)}}
{c_{1}^{(j)}\psi_{k-1}^{(j)}(z)}\,,\qquad 1\leq j \leq m-1\,.
\end{equation}
\end{defn}
Notice that $\varphi_{k-1}^{(1)}=\varphi_{k-1}$, where
$\varphi_{k-1}$ was previously defined in (\ref{eqdefvarphi}).
\begin{theo}\label{relasympstf}
Let $S = {\mathcal{N}}^*(\sigma_1,\ldots,\sigma_m)$ and $\Lambda
\subset {\mathbb{Z}}_+^{m}(\circledast;p_1,\ldots,p_m)$ be a
sequence of distinct multi-indices such that for all ${\bf n} \in
\Lambda, n_{1} - n_{m} \leq  C,$ where $C$ is a constant. Then,
for each $k\in\{0,1,\ldots,m\}$,
\begin{equation}\label{asintRPsi}
\lim_{\mathbf{n}\in\Lambda}\frac{\widetilde{\Psi}_{\mathbf{n},k}(z)}
{\Psi_{\mathbf{n},k}(z)}=G_{k}(z;p_{1},\ldots,p_{m})\,, \quad K
\subset \overline{\mathbb{C}} \setminus
(\supp{\sigma_{k}}\cup\supp{\sigma_{k+1}})\,,
\end{equation}
where $G_{k}$ is analytic and never vanishes in the indicated
region. For each $k = \{0,\ldots,m-1\}$ and all sufficiently large
$|\mathbf{n}|, \mathbf{n} \in \Lambda$, $ \widetilde{\Psi}_{{\bf
n},k}$ has exactly $N_{{\bf n},k+1}= n_{k+1}+\cdots + n_m$ zeros
in $\overline{\mathbb{C}} \setminus \supp{\sigma_{k}}$,
$\supp{\sigma_{k+1}}$ is an attractor of the zeros of
$\{\widetilde{\Psi}_{{\bf n},k}\}, {\bf n} \in \Lambda$, in this
region, and each point in $\supp{\sigma_{k+1}} \setminus
\widetilde{\Delta}_{k+1}$ is a $1$ attraction point of zeros of
$\{\widetilde{\Psi}_{{\bf n},k}\}, {\bf n} \in \Lambda$. When the
coefficients of the polynomials $p_k, k=1,\ldots,m,$ are real, all
the statements above  remain valid for $\Lambda \subset
\mathbb{Z}_+^m(\circledast).$ An expression for $G_{k}$ is given
in $(\ref{eqcasokigual2})$-$(\ref{eqcasokmayortres})$ below.
\end{theo}
{\bf Proof.} For $k=0$, $(\ref{asintRPsi})$  is (\ref{eq:lim})
since $\widetilde{\Psi}_{\mathbf{n},0}=\widetilde{Q}_{\mathbf{n}}$
and $\Psi_{\mathbf{n},0}=Q_{\mathbf{n}}$; therefore,
\[
G_{0}(z;p_{1},\ldots,p_{m})={\mathcal{F}}(z;p_1,\ldots,p_m)\,.
\]

By (\ref{eq:32}), we know that
\[
\lim_{\mathbf{n}\in\Lambda}\frac{\lambda_{\mathbf{n}}^*
\varepsilon_{\mathbf{n}_0,k-1}K_{\mathbf{n}_0,k-1}^2
(Q_{\mathbf{n}_0,k-1}R_{\mathbf{n},k-1})(z)}{Q_{\mathbf{n}_0,k}(z)}
=
\]
\[
\left\{
\begin{array}{ll}
\frac{1}{\sqrt{(z-b_{1})(z-a_{1})}} p_{\Lambda}(\varphi_{1}(z))\,,& k=2\,, \\
\frac{1}{\sqrt{(z-b_{k-1})(z-a_{k-1})}}p_{\Lambda}(\delta\varphi_{k-1}(z))\,, & k=3,\ldots,m\,, \\
\end{array} \right.
\]
uniformly on compact subsets of
$\mathbb{C}\setminus(\supp{\sigma_{k-1}}\cup\supp{\sigma_{k}})$.
Also, see (\ref{eq:x}),
\[
\lim_{\mathbf{n}\in\Lambda}\varepsilon_{\mathbf{n}_{0},k-1}
h_{\mathbf{n}_{0},k}(z)=\frac{1}{\sqrt{(z-b_{k-1})(z-a_{k-1})}}\,,
\qquad K \subset
\overline{\mathbb{C}}\setminus\supp{\sigma_{k-1}}.
\]
Thus, since
$\lim_{\mathbf{n}\in\Lambda}\lambda_{\mathbf{n}}^{*}=c$, we
conclude that
\begin{equation}\label{parte1}
\lim_{\mathbf{n}\in\Lambda}\frac{R_{\mathbf{n},k-1}(z)
}{\Psi_{\mathbf{n}_{0},k-1}(z)}
=\lim_{\mathbf{n}\in\Lambda}K_{\mathbf{n}_0,k-1}^2\frac{
(Q_{\mathbf{n}_0,k-1}R_{\mathbf{n},k-1})(z)}
{(h_{\mathbf{n}_{0},k}Q_{\mathbf{n}_0,k})(z)}=
\end{equation}
\[
\left\{
\begin{array}{ll}
p_{\Lambda}(\varphi_{1}(z))/c\,,& k=2\,, \\
p_{\Lambda}(\delta\varphi_{k-1}(z))/c\,, & k=3,\ldots,m\,, \\
\end{array} \right.
\]
uniformly on compact subsets of
$\mathbb{C}\setminus(\supp{\sigma_{k-1}}\cup\supp{\sigma_{k}})$.

Recall that $\mathbf{n}_{j}= (n_1 - \deg(p_2\cdots p_m) +j, n_2 -
\deg (p_3\cdots p_m),\ldots, n_m)$. It is easy to see that
\[
\frac{\Psi_{\mathbf{n}_{0},k-1}}{\Psi_{\mathbf{n}_{j},k-1}}=
\frac{Q_{\mathbf{n}_0,k}}{Q_{\mathbf{n}_{j},k}}
\frac{Q_{\mathbf{n}_{j},k-1}}{Q_{\mathbf{n}_0,k-1}}
\frac{\varepsilon_{\mathbf{n}_{0},k-1}
h_{\mathbf{n}_{0},k}}{\varepsilon_{\mathbf{n}_{j},k-1}
h_{\mathbf{n}_{j},k}}
\frac{\varepsilon_{\mathbf{n}_{j},k-1}}{\varepsilon_{\mathbf{n}_{0},k-1}}
\frac{K_{\mathbf{n}_{j},k-1}^2}{K_{\mathbf{n}_{0},k-1}^2}\,.
\]
From this expression, applying Proposition \ref{cor3} and
(\ref{ratioepsilons}), we obtain that the following limit holds
uniformly on compact subsets of
$\mathbb{C}\setminus(\supp{\sigma_{k-1}}\cup\supp{\sigma_{k}})$
\[
\lim_{\mathbf{n}\in\Lambda}
\frac{\Psi_{\mathbf{n}_{0},k-1}(z)}{\Psi_{\mathbf{n}_{j},k-1}(z)}=
(\Delta_{k-1,1}\cdots\Delta_{1,1})^{j}
\Big(\frac{\widetilde{F}_{k-1}^{(1)}(z)}{\widetilde{F}_{k}^{(1)}(z)}\Big)^{j}
(\kappa^{(1)}_{1}\cdots\kappa^{(1)}_{k-1})^{2j}\,.
\]
Now, from (\ref{eq:xc}) and (\ref{eq:F}), we have
\[\frac{\widetilde{F}_{k-1}^{(1)}(z)}{\widetilde{F}_{k}^{(1)}(z)}
=\frac{c_{k}^{(1)}}{c_{k-1}^{(1)}}\sg{\psi_{k-1}^{(1)}(\infty)}\psi_{k-1}^{(1)}(z)\,,\]
and from (\ref{eq:xc})
\[
(\kappa^{(1)}_{1}\cdots\kappa^{(1)}_{k-1})^{2}=c_{1}^{(1)}\frac{c_{k-1}^{(1)}}{c_{k}^{(1)}}\,.
\]
Thus,
\[
\lim_{\mathbf{n}\in\Lambda}
\frac{\Psi_{\mathbf{n}_{0},k-1}(z)}{\Psi_{\mathbf{n}_{j},k-1}(z)}=(\Delta_{k-1,1}\cdots\Delta_{1,1})^{j}
(c_{1}^{(1)}\sg{\psi_{k-1}^{(1)}(\infty)}\psi_{k-1}^{(1)}(z))^{j}\,.
\]

Set
\begin{equation}\label{definitionxi}
\Xi_{k}:=(\Delta_{k-1,1}\cdots\Delta_{1,1})^{\deg(p_2\cdots
p_{m})} \cdots
(\Delta_{k-1,m-1}\cdots\Delta_{1,m-1})^{\deg(p_{m})}\,.
\end{equation}
Using the same arguments above, on an appropriate consecutive
collection of multi-indices, one proves that
\[
\lim_{\mathbf{n}\in\Lambda}
\frac{\Psi_{\mathbf{n}_{0},k-1}(z)}{\Psi_{\mathbf{n},k-1}(z)} =
\Xi_{k}
\prod_{j=1}^{m-1}\frac{1}{(\varphi_{k-1}^{(j)}(z))^{\deg(p_{j+1}\cdots
p_{m})}}\,,
\]
uniformly on compact subsets of
$\mathbb{C}\setminus(\supp{\sigma_{k-1}}\cup\supp{\sigma_{k}})$.
Therefore, writing
\[\frac{R_{\mathbf{n},k-1}(z)}
{\Psi_{\mathbf{n},k-1}(z)}=\frac{R_{\mathbf{n},k-1}(z)}
{\Psi_{\mathbf{n}_{0},k-1}(z)}
\frac{\Psi_{\mathbf{n}_{0},k-1}(z)}{\Psi_{\mathbf{n},k-1}(z)}\,,
\]
using the expression of $p_{\Lambda}$, applying (\ref{parte1}) and
Lemma \ref{lemmarelation}, for $k=2$ we get
\[
\lim_{\mathbf{n}\in\Lambda}\frac{\widetilde{\Psi}_{\mathbf{n},1}(z)}
{\Psi_{\mathbf{n},1}(z)}= \Xi_{2}\prod_{\nu=1}^{l_{1}}
(\varphi_{1}(z)-\varphi_{0}(z_{1,\nu}))^{\tau_{1,\nu}}
\prod_{\nu=1}^{l_{2}}\Big(\frac{1}{\varphi_{1}(z)}
\frac{\varphi_{1}(z)-\varphi_{1}(z_{2,\nu})}{z-z_{2,\nu}}
\Big)^{\tau_{2,\nu}}
\]
\begin{equation}\label{eqcasokigual2}
\times\prod_{j=3}^{m}\prod_{\nu=1}^{l_{j}}\Big(
\frac{\varphi_{1}(z)-\delta\varphi_{j-1}(z_{j,\nu})}
{\varphi_{1}^{(j-1)}(z)}\Big)^{\tau_{j,\nu}}
\end{equation}
uniformly on compact subsets of
$\overline{\mathbb{C}}\setminus(\supp{\sigma_{1}}\cup\supp{\sigma_{2}})$,
and for $k\geq 3$ we obtain
\[
\lim_{\mathbf{n}\in\Lambda}\frac{\widetilde{\Psi}_{\mathbf{n},k-1}(z)}
{\Psi_{\mathbf{n},k-1}(z)}= \Xi_{k}\prod_{\nu=1}^{l_{1}}
(\delta\varphi_{k-1}(z)-\varphi_{0}(z_{1,\nu}))^{\tau_{1,\nu}}
\prod_{\nu=1}^{l_{2}}\Big(
\frac{\delta\varphi_{k-1}(z)-\varphi_{1}(z_{2,\nu})}{\varphi_{k-1}(z)}
\Big)^{\tau_{2,\nu}}
\]
\begin{equation}\label{eqcasokmayortres}
\times\prod_{\nu=1}^{l_{k}}\Big(
\frac{\delta\varphi_{k-1}(z)-\delta\varphi_{k-1}(z_{k,\nu})}
{\varphi_{k-1}^{(k-1)}(z)(z-z_{k,\nu})}\Big)^{\tau_{k,\nu}}\prod_{j=3,
j\neq k}^{m}\prod_{\nu=1}^{l_{j}}\Big(
\frac{\delta\varphi_{k-1}(z)-\delta\varphi_{j-1}(z_{j,\nu})}
{\varphi_{k-1}^{(j-1)}(z)}\Big)^{\tau_{j,\nu}}
\end{equation}
uniformly on compact subsets of
$\overline{\mathbb{C}}\setminus(\supp{\sigma_{k-1}}\cup\supp{\sigma_{k}})$.
Therefore, (\ref{asintRPsi}) is proved.

From the expression of the limit functions one sees that $G_{k}$
does not vanish in $\overline{\mathbb{C}} \setminus
(\supp{\sigma_{k}}\cup\supp{\sigma_{k+1}})$. The statements
concerning the number of zeros of $\widetilde{\Psi}_{{\bf n},k}$
for $k\in \{0,\ldots,m-1\}$ and their limit behavior follows at
once from (\ref{asintRPsi}), on account of the argument principle
and the corresponding behavior of the zeros of the polynomials
$Q_{{\bf n},k+1}$ described in Proposition \ref{pro1}. Recall that
the zeros of $Q_{{\bf n},k+1}$ are those of $\Psi_{{\bf n},k}$ in
$\overline{\mathbb{C}} \setminus \supp{\sigma_{k}}$.

Now, let us assume that the coefficients of the polynomials
$p_{k}$ are real and $\Lambda \subset\mathbb{Z}_+^m(\circledast)$.
Since
\[
\frac{\Psi_{\mathbf{n}^{l},k-1}}{\Psi_{\mathbf{n},k-1}}=
\frac{Q_{\mathbf{n}^{l},k}}{Q_{\mathbf{n},k}}
\frac{Q_{\mathbf{n},k-1}}{Q_{\mathbf{n}^{l},k-1}}
\frac{\varepsilon_{\mathbf{n}^{l},k-1}
h_{\mathbf{n}^{l},k}}{\varepsilon_{\mathbf{n},k-1}
h_{\mathbf{n},k}}
\frac{\varepsilon_{\mathbf{n},k-1}}{\varepsilon_{\mathbf{n}^{l},k-1}}
\frac{K_{\mathbf{n},k-1}^2}{K_{\mathbf{n}^{l},k-1}^2}\,,
\]
applying (\ref{eq:xk}), (\ref{eq:xb}), (\ref{eq:x}), and
(\ref{ratioepsilons}), we conclude that the ratio asymptotic
\[
\lim_{{\bf n} \in
\Lambda}\frac{\Psi_{\mathbf{n}^{l},k-1}(z)}{\Psi_{\mathbf{n},k-1}(z)}
\,,\qquad K \subset \mathbb{C} \setminus
(\supp{\sigma_{k-1}}\cup\supp{\sigma_{k}})\,,
\]
holds and the limit does not vanish in the indicated region.

Since each measure $p_{k}\,\sigma_{k}$ is real with constant sign,
we can define the polynomials $\widetilde{Q}_{{\bf n},k}, 1\leq
k\leq m,$ as the monic polynomials of degree $N_{{\bf n},k}$ whose
simple zeros are located at the points where
$\widetilde{\Psi}_{n,k-1}$ vanishes on $\Delta_{k}$. Let
$\widetilde{Q}_{{\bf n},0}\equiv\widetilde{Q}_{{\bf n},m+1}\equiv
1$. We also introduce the associated notions
\begin{equation}\label{definhtilde}
\widetilde{H}_{{\bf n},k}:=  \frac{\widetilde{Q}_{{\bf
n},k-1}\widetilde{\Psi}_{{\bf n},k-1}}{\widetilde{Q}_{{\bf
n},k}}\,, \qquad k=1,\ldots,m+1\,,
\end{equation}
$\widetilde{\varepsilon}_{\mathbf{n},k}$ as the sign of $
{\widetilde{H}_{{\bf
n},k}(x)p_{k}(x)d\sigma_k(x)}/{\widetilde{Q}_{{\bf n},k-1}(x)
\widetilde{Q}_{{\bf n},k+1}(x)}$ on $\supp{\sigma_k}$, and
\begin{equation}\label{defnKnktilde}
\widetilde{K}_{{\bf n},k}:= \left( \int \widetilde{Q}_{{\bf
n},k}^2(x)
\frac{\widetilde{\varepsilon}_{\mathbf{n},k}\widetilde{H}_{\mathbf{n},k}(x)p_{k}(x)d\sigma_k(x)}{\widetilde{Q}_{{\bf
n},k-1}(x)\widetilde{Q}_{{\bf n},k+1}(x)} \right)^{-1/2}\;.
\end{equation}
The formulas (\ref{eq:xk}), (\ref{eq:xb}), (\ref{eq:x}) and
(\ref{ratioepsilons}) are independent of the orthogonality
measures, hence
\[
\lim_{{\bf n} \in
\Lambda}\frac{\widetilde{\Psi}_{\mathbf{n}^{l},k-1}(z)}
{\widetilde{\Psi}_{\mathbf{n},k-1}(z)}=\lim_{{\bf n} \in
\Lambda}\frac{\Psi_{\mathbf{n}^{l},k-1}(z)}
{{\Psi}_{\mathbf{n},k-1}(z)}\,.
\]
Applying the same argument used in the last two paragraphs of the
proof of Theorem \ref{teo:1}, we conclude that (\ref{asintRPsi})
is valid for $\Lambda\subset\mathbb{Z}_+^m(\circledast)$.

The rest of the statements regarding the zeros of $
\widetilde{\Psi}_{{\bf n},k}$   and their limit behavior follows
as in the case of polynomials with complex coefficients. \hfill
$\Box$
\begin{co} \label{cor5.1}
Let $S = {\mathcal{N}}^*(\sigma_1,\ldots,\sigma_m)$. Consider the
perturbed Nikishin system
$\mathcal{N}(\frac{p_1}{q_1}\sigma_1,\ldots,\frac{p_m}{q_m}\sigma_m)$,
where $p_k,q_k$ denote relatively prime polynomials whose zeros
lie in  $\mathbb{C} \setminus \cup_{k=1}^m \Delta_k$. Let $\Lambda
\subset {\mathbb{Z}}^{m}_+(\circledast;p_1q_1,\ldots,p_mq_m)$ be a
sequence of distinct multi-indices such that for all ${\bf n} \in
\Lambda$, $n_{1} - n_{m} \leq C,$ where $C$ is a constant. Let
$\widetilde{Q}_{{\bf n}}$ be the monic multiple orthogonal
polynomial of smallest degree relative to the Nikishin system
$\mathcal{N}(\frac{p_1}{q_1}\sigma_1,\ldots,\frac{p_m}{q_m}\sigma_m)$
and ${\bf n}$, whereas $\widetilde{\Psi}_{{\bf n},k}, 0\leq k\leq
m,$ denote the second type functions defined in
$(\ref{defnPsinktilde})$, with $p_{k}$ replaced by $p_{k}/q_{k}$.
Then, for each $k\in\{0,\ldots,m\}$,
\begin{equation}\label{eq:relativesecondtyperational}
\lim_{{\bf n} \in
\Lambda}\frac{\widetilde{\Psi}_{n,k}(z)}{\Psi_{n,k}(z)}
=\frac{G_{k}(z;p_{1},\ldots,p_{m})}
{G_{k}(z;q_{1},\ldots,q_{m})}\,,\qquad
K\subset\overline{\mathbb{C}} \setminus
(\supp{\sigma_{k}}\cup\supp{\sigma_{k+1}})\,.
\end{equation}
For each $k = \{0,\ldots,m-1\}$ and all sufficiently large
$|\mathbf{n}|, \mathbf{n} \in \Lambda$, $ \widetilde{\Psi}_{{\bf
n},k}$ has exactly $N_{{\bf n},k+1}$ zeros in
$\overline{\mathbb{C}} \setminus \supp{\sigma_{k}}$,
$\supp{\sigma_{k+1}}$ is an attractor of the zeros of
$\{\widetilde{\Psi}_{{\bf n},k}\}, {\bf n} \in \Lambda$, in this
region, and each point in $\supp{\sigma_{k+1}} \setminus
\widetilde{\Delta}_{k+1}$ is a $1$ attraction point of zeros of
$\{\widetilde{\Psi}_{{\bf n},k}\}, {\bf n} \in \Lambda$. When the
polynomials $p_{k}, q_{k}, k=1,\ldots,m$, have real coefficients,
all the statements remain valid when
$\Lambda\subset{\mathbb{Z}}^{m}_+(\circledast)$.
\end{co}
{\bf Proof.} We consider the auxiliary Nikishin system
\[
S_{1}:= {\mathcal{N}}\Big(\frac{\sigma_1}{|q_{1}|^{2}},\ldots,
\frac{\sigma_m}{|q_{m}|^{2}}\Big)\,,
\]
and define the related second type functions
\[\Psi_{n,0}^{*}(z):=Q_{{\bf n}}^{*}(z)\,,\]
\[
\Psi_{n,k}^{*}(z):=\int\frac{\Psi_{n,k-1}^{*}(x)}{z-x}
\frac{d\sigma_{k}(x)}{|q_{k}(x)|^{2}}\,,\quad 1\leq k\leq m\,,
\]
where $Q_{{\bf n}}^{*}$ denotes the multiple orthogonal polynomial
associated to $S_{1}$ and ${\bf n}$.

Notice that if we perturb the generator of system $S_{1}$
multiplying the $k$-th measure by the real polynomial
$|q_{k}|^{2}$ we get the generator of the original Nikishin system
$S$. Thus, applying Theorem \ref{relasympstf}, we obtain that for
all $k\in\{0,\ldots,m\}$
\[
\lim_{{\bf n} \in \Lambda}\frac{\Psi_{n,k}(z)}{\Psi_{n,k}^{*}(z)}
=G_{k}(z;|q_{1}|^{2},\ldots,|q_{m}|^{2})\,, \quad K \subset
\overline{\mathbb{C}} \setminus
(\supp{\sigma_{k}}\cup\supp{\sigma_{k+1}})\,.
\]
The perturbed system $S_{2}:=
\mathcal{N}(\frac{p_1}{q_1}\sigma_1,\ldots,\frac{p_m}{q_m}\sigma_{m})$
can be written as
\[
S_{2}={\mathcal{N}}\Big(p_{1}\overline{q}_1
\frac{\sigma_1}{|q_{1}|^{2}},\ldots,
p_{m}\overline{q}_{m}\frac{\sigma_{m}}{|q_{m}|^{2}}\Big)\,.
\]
Therefore, employing the same argument
\[
\lim_{{\bf n} \in
\Lambda}\frac{\widetilde{\Psi}_{n,k}(z)}{\Psi_{n,k}^{*}(z)}
=G_{k}(z;p_{1}\overline{q}_1,\ldots,p_{m}\overline{q}_{m})\,,
\quad K \subset \overline{\mathbb{C}} \setminus
(\supp{\sigma_{k}}\cup\supp{\sigma_{k+1}})\,.
\]

We conclude that
\[
\lim_{{\bf n} \in
\Lambda}\frac{\widetilde{\Psi}_{n,k}(z)}{\Psi_{n,k}(z)}
=\frac{G_{k}(z;p_{1}\overline{q}_1,\ldots,p_{m}\overline{q}_{m})}
{G_{k}(z;q_{1}\overline{q}_1,\ldots,q_{m}\overline{q}_{m})}=
\frac{G_{k}(z;p_{1},\ldots,p_{m})}
{G_{k}(z;q_{1},\ldots,q_{m})}\,,
\]
uniformly on compact subsets of $\overline{\mathbb{C}} \setminus
(\supp{\sigma_{k}}\cup\supp{\sigma_{k+1}})$. The statements
concerning the zeros can be proved  as in the case of polynomial
perturbation.

When the polynomials $p_{k}, q_{k}, k=1,\ldots,m$, have real
coefficients, it follows from Theorem \ref{relasympstf} that
(\ref{eq:relativesecondtyperational}) remains valid for
$\Lambda\subset{\mathbb{Z}}^{m}_+(\circledast)$. The statements
concerning the zeros are derived immediately. \hfill $\Box$

\section{Relative asymptotics for the polynomials $Q_{{\bf n},k}$}
In this section, we will restrict our attention to the case when
the polynomials $p_k,q_k, k=1,\ldots,m,$ have real coefficients
(and of course their zeros lie in $\mathbb{C}\setminus
\cup_{k=1}^m \Delta_{k}$). Accordingly, we use the objects
$\widetilde{Q}_{{\bf n},k}, \widetilde{H}_{{\bf n},k},
\widetilde{K}_{{\bf n},k},$ and $\widetilde{\varepsilon}_{{\bf
n},k},$ introduced at the end of the proof of Theorem
\ref{relasympstf} (see (\ref{definhtilde}) and
(\ref{defnKnktilde})). Here, we study the relative asymptotic of
the ratios $\widetilde{Q}_{{\bf n},k}/Q_{{\bf n},k}$.

\begin{lem}\label{relationepsilons}
For any $\mathbf{n}\in\mathbb{Z}_{+}^{m}(\circledast)$
\begin{equation}\label{eqsignos}
\frac{\varepsilon_{\mathbf{n},k}}{\widetilde{\varepsilon}_{\mathbf{n},k}}=\prod_{i=1}^{k}\sign{p_{i},\supp{\sigma_{i}}}\,.
\end{equation}
\end{lem}
{\bf Proof.} By definition $\varepsilon_{\mathbf{n},k}$ is the
sign of $H_{{\bf n},k}(x)d\sigma_k(x)/Q_{{\bf n},k-1}(x)Q_{{\bf
n},k+1}(x)$ on $\supp{\sigma_k}$ and
$\widetilde{\varepsilon}_{\mathbf{n},k}$ is the sign of
$\widetilde{H}_{{\bf
n},k}(x)p_{k}(x)d\sigma_k(x)/\widetilde{Q}_{{\bf n},k-1}(x)
\widetilde{Q}_{{\bf n},k+1}(x)$ on $\supp{\sigma_k}$. If $k=1$
these measures reduce respectively to $d\sigma_1(x)/Q_{{\bf
n},2}(x)$ and $p_1(x)d\sigma_1(x)/\widetilde{Q}_{{\bf n},2}(x)$.
Since $Q_{{\bf n},2}$ and $\widetilde{Q}_{{\bf n},2}$ are monic
polynomials of the same degree and their zeros are located in
$\Delta_{2}$, which is disjoint with $\supp{\sigma_1}$, it follows
that $Q_{{\bf n},2}$ and $\widetilde{Q}_{{\bf n},2}$ have the same
sign on $\supp{\sigma_1}$. Therefore,
\[
\frac{\varepsilon_{\mathbf{n},1}}{\widetilde{\varepsilon}_{\mathbf{n},1}}=\sign{p_{1},\supp{\sigma_{1}}}\,.
\]
To conclude the proof we show that
\[
\frac{\varepsilon_{\mathbf{n},k}}{\widetilde{\varepsilon}_{\mathbf{n},k}}=\sign{p_{k},\supp{\sigma_{k}}}
\frac{\varepsilon_{\mathbf{n},k-1}}{\widetilde{\varepsilon}_{\mathbf{n},k-1}}\,.
\]

Notice that $Q_{{\bf n},k-1}$ and $\widetilde{Q}_{{\bf n},k-1}$
have the same sign on $\supp{\sigma_{k}}$ by an argument similar
to the one explained above. The same holds for $Q_{{\bf n},k+1}$
and $\widetilde{Q}_{{\bf n},k+1}$. Therefore
\[
\frac{\varepsilon_{\mathbf{n},k}}{\widetilde{\varepsilon}_{\mathbf{n},k}}=\frac{\sign{H_{{\bf
n},k},\supp{\sigma_{k}}}}{\sign{p_{k}\widetilde{H}_{{\bf
n},k},\supp{\sigma_{k}}}}\,.
\]

By (\ref{eq:18}), we know  that
\[
H_{\mathbf{n},k}(x)=\int_{\Delta_{k-1}}\frac{Q_{\mathbf{n},k-1}^{2}(t)}{x-t}\frac{H_{\mathbf{n},k-1}(t)
d\sigma_{k-1}(t)}{Q_{\mathbf{n},k-2}(t)Q_{\mathbf{n},k}(t)}\,,
\]
and
\[
\widetilde{H}_{\mathbf{n},k}(x)=\int_{\Delta_{k-1}}\frac{\widetilde{Q}_{\mathbf{n},k-1}^{2}(t)}{x-t}\frac{\widetilde{H}_{\mathbf{n},k-1}(t)
p_{k-1}(t)d\sigma_{k-1}(t)}{\widetilde{Q}_{\mathbf{n},k-2}(t)\widetilde{Q}_{\mathbf{n},k}(t)}\,.
\]
Consequently,
\[
\frac{\sign{H_{{\bf
n},k},\supp{\sigma_{k}}}}{\sign{\widetilde{H}_{{\bf
n},k},\supp{\sigma_{k}}}}=\frac{\varepsilon_{\mathbf{n},k-1}}{\widetilde{\varepsilon}_{\mathbf{n},k-1}}\,,
\]
and the claim follows.\hfill $\Box$

We are ready to state and prove
\begin{theo}\label{theorelativepol}
Let $S = {\mathcal{N}}^*(\sigma_1,\ldots,\sigma_m)$ and $\Lambda
\subset {\mathbb{Z}}_+^{m}(\circledast)$ be a sequence of distinct
multi-indices such that for all ${\bf n} \in \Lambda, n_{1} -
n_{m} \leq  C,$ where $C$ is a constant. Assume that the
polynomials $p_k, k=1,\ldots,m,$ have real coefficients. For each
$k\in\{1,\ldots,m\}$,
\begin{equation}\label{asintrelativtilde}
\lim_{\mathbf{n}\in\Lambda}\frac{\widetilde{Q}_{\mathbf{n},k}(z)}
{Q_{\mathbf{n},k}(z)}=\mathcal{F}_{k}(z; p_1 ,\ldots,p_m)\,, \quad
K \subset \overline{\mathbb{C}} \setminus \supp{\sigma_{k}}\,,
\end{equation}
where $\mathcal{F}_{k}(z; p_1 ,\ldots,p_m)$ is analytic and never
vanishes in $\overline{\mathbb{C}} \setminus \supp{\sigma_{k}}$
and
\begin{equation}\label{asintrelativKcuad}
\lim_{\mathbf{n}\in\Lambda}\frac{\widetilde{K}_{\mathbf{n},k}^{2}}
{K_{{\bf n},k}^{2}}=\frac{\prod_{i=1}^{k}
\sign{p_{i},\supp{\sigma_{i}}}}{G_{k}(\infty;p_{1},\ldots,p_{m})}.
\end{equation}
For  $k\in\{1,\ldots,m-1\}$ and
$z\in\overline{\mathbb{C}}\setminus
(\supp{\sigma_{k}}\cup\supp{\sigma_{k+1}})$
\begin{equation}\label{relationFk}
\mathcal{F}_{k+1}(z;p_{1},\ldots,p_{m})=
\prod_{i=0}^{k}\frac{G_{i}(z;p_{1},\ldots,p_{m})}{G_{i}(\infty;p_{1},\ldots,p_{m})}\,,
\end{equation}
where $G_{i}(z;p_{1},\ldots,p_{m})$ is the function given in
$(\ref{asintRPsi})$.
\end{theo}

{\bf Proof.} If $\Lambda \subset
\mathbb{Z}_+^m(\circledast;p_1,\ldots,p_m)$, from (\ref{eq:32})
and Lemma \ref{lemmarelation}, we have that
\[ \lim_{\mathbf{n} \in \Lambda} \lambda_{\mathbf{n}}^*
\varepsilon_{\mathbf{n}_0,k-1}K_{\mathbf{n}_0,k-1}^2  \frac{
Q_{\mathbf{n}_0,k-1}(z)(p_{k}\cdots
p_{m})(z)\widetilde{\Psi}_{{\bf n},k-1}(z)}{Q_{\mathbf{n}_0,k}(z)}
=
\]
\begin{equation}\label{otraparapsitilde}
\left\{
\begin{array}{ll}
\frac{1}{\sqrt{(z-b_{1})(z-a_{1})}} p_{\Lambda}(\varphi_{1}(z))\,,& k=2\,, \\
\frac{1}{\sqrt{(z-b_{k-1})(z-a_{k-1})}}p_{\Lambda}(\delta\varphi_{k-1}(z))\,, & k=3,\ldots,m\,, \\
\end{array} \right.
\end{equation}
By Proposition \ref{pro1}, we know that
\begin{equation}\label{eqasinttilde}
\lim_{{\bf n}\in
{\Lambda}}\widetilde{\varepsilon}_{\mathbf{n},k}\widetilde{K}_{{\bf
n},k}^{2}\widetilde{H}_{\mathbf{n},k+1}(z) = \frac{1}{\sqrt{(z -
b_{k})(z- a_{k})}}\,, \qquad K \subset \overline{\mathbb{C}}
\setminus \supp{\sigma_k}\,,
\end{equation}
where $[a_{k},b_{k}]=\widetilde{\Delta}_k$. Formula
(\ref{definhtilde}) implies
\[
\frac{\lambda_{\mathbf{n}}^*
\varepsilon_{\mathbf{n}_0,k-1}K_{\mathbf{n}_0,k-1}^2
Q_{\mathbf{n}_0,k-1}(z)(p_{k}\cdots
p_{m})(z)\widetilde{\Psi}_{{\bf
n},k-1}(z)}{\widetilde{\varepsilon}_{\mathbf{n},k-1}\widetilde{K}_{{\bf
n},k-1}^{2}\widetilde{H}_{\mathbf{n},k}(z)Q_{\mathbf{n}_0,k}(z)}=
\]
\begin{equation}\label{eq4}
\lambda_{\mathbf{n}}^* \frac{\varepsilon_{\mathbf{n}_0,k-1}}
{\widetilde{\varepsilon}_{\mathbf{n},k-1}}
\frac{K_{\mathbf{n}_0,k-1}^2}{\widetilde{K}_{{\bf
n},k-1}^{2}}\frac{Q_{\mathbf{n}_0,k-1}(z)}
{\widetilde{Q}_{\mathbf{n},k-1}(z)}\frac{\widetilde{Q}_{\mathbf{n},k}(z)}
{Q_{\mathbf{n}_{0},k}(z)}(p_{k}\cdots p_{m})(z)\,.
\end{equation}

Using (\ref{otraparapsitilde}), (\ref{eqasinttilde}), and
(\ref{eq4}), we obtain
\[
\lim_{{\bf n}\in {\Lambda}}\lambda_{\mathbf{n}}^*
\frac{\varepsilon_{\mathbf{n}_0,k-1}}
{\widetilde{\varepsilon}_{\mathbf{n},k-1}}
\frac{K_{\mathbf{n}_0,k-1}^2}{\widetilde{K}_{{\bf
n},k-1}^{2}}\frac{Q_{\mathbf{n}_0,k-1}(z)}
{\widetilde{Q}_{\mathbf{n},k-1}(z)}\frac{\widetilde{Q}_{\mathbf{n},k}(z)}
{Q_{\mathbf{n}_{0},k}(z)}(p_{k}\cdots p_{m})(z) =
\]
\begin{equation}\label{eqasintotica}
\left\{
\begin{array}{ll}
p_{\Lambda}(\varphi_{1}(z))\,,& k=2\,, \\
p_{\Lambda}(\delta\varphi_{k-1}(z))\,, & k=3,\ldots,m\,, \\
\end{array} \right.
\end{equation}
Using the results on ratio asymptotic for the constants $K_{{\bf
n},k}, \widetilde{K}_{{\bf n},k}$ and the polynomials
$Q_{\mathbf{n},k}, \widetilde{Q}_{\mathbf{n},k}$, it follows that
(\ref{eqasintotica}) is also valid for
$\Lambda\subset{\mathbb{Z}}^{m}_+(\circledast)$.

Since
\[\frac{\varepsilon_{\mathbf{n}_{0},k-1}}{\widetilde{\varepsilon}_{\mathbf{n},k-1}}=
\frac{\varepsilon_{\mathbf{n}_{0},k-1}}{{\varepsilon}_{\mathbf{n},k-1}}
\frac{\varepsilon_{\mathbf{n},k-1}}{\widetilde{\varepsilon}_{\mathbf{n},k-1}}\,,\]
applying Lemma \ref{relationepsilons}, (\ref{ratioepsilons}), and
(\ref{definitionxi}), we obtain
\begin{equation}\label{ratioepsilons2}
\frac{\varepsilon_{\mathbf{n}_{0},k-1}}{\widetilde{\varepsilon}_{\mathbf{n},k-1}}=
\Xi_{k}\prod_{i=1}^{k-1}\sign{p_{i},\supp{\sigma_{i}}}\,.
\end{equation}

We have
\begin{equation}\label{ratioks}
\frac{K_{\mathbf{n}_0,k-1}^2}{\widetilde{K}_{{\bf n},k-1}^{2}}=
\frac{K_{\mathbf{n}_0,k-1}^2}{K_{\mathbf{n},k-1}^2}\frac{K_{\mathbf{n},k-1}^2}{\widetilde{K}_{{\bf
n},k-1}^{2}}\,,
\end{equation}
and by (\ref{eq:xk})
\begin{equation}\label{limitratioks}
\lim_{\mathbf{n}\in\Lambda}\frac{K_{\mathbf{n}_0,k-1}^2}{K_{\mathbf{n},k-1}^2}=
\prod_{i=1}^{m-1}(\kappa^{(i)}_{1}\cdots\kappa^{(i)}_{k-1})^{-2\deg(p_{i+1}\cdots
p_{m})}\,.
\end{equation}

Write
\begin{equation}\label{ratioq}
\frac{Q_{\mathbf{n}_0,k-1}(z)}{\widetilde{Q}_{\mathbf{n},k-1}(z)}=
\frac{Q_{\mathbf{n}_0,k-1}(z)}{Q_{\mathbf{n},k-1}(z)}
\frac{Q_{\mathbf{n},k-1}(z)}{\widetilde{Q}_{\mathbf{n},k-1}(z)}\,,
\end{equation}
and
\begin{equation}\label{ratioq2}
\frac{\widetilde{Q}_{\mathbf{n},k}(z)}{Q_{\mathbf{n}_{0},k}(z)}=
\frac{\widetilde{Q}_{\mathbf{n},k}(z)}{Q_{\mathbf{n},k}(z)}\frac{Q_{\mathbf{n},k}(z)}{Q_{\mathbf{n}_{0},k}(z)}\,.
\end{equation}
Notice that
\begin{equation}\label{limitratioq}
\lim_{\mathbf{n}\in\Lambda}\frac{Q_{\mathbf{n}_0,k-1}(z)}{Q_{\mathbf{n},k-1}(z)}=
\prod_{i=1}^{m-1}(\widetilde{F}_{k-1}^{(i)}(z))^{-deg(p_{i+1}\cdots
p_{m})}\,.
\end{equation}
\begin{equation}\label{limitratioq2}
\lim_{\mathbf{n}\in\Lambda}\frac{Q_{\mathbf{n},k}(z)}{Q_{\mathbf{n}_{0},k}(z)}=
\prod_{i=1}^{m-1}(\widetilde{F}_{k}^{(i)}(z))^{deg(p_{i+1}\cdots
p_{m})}\,.
\end{equation}

From (\ref{eq:F}) and (\ref{eq:xc}) it follows that
\[
\frac{\widetilde{F}_{k}^{(i)}(z)}{\widetilde{F}_{k-1}^{(i)}(z)}=
\frac{c_{k-1}^{(i)}}{c_{k}^{(i)}}\frac{\sg{\psi_{k-1}^{(i)}(\infty)}}{\psi_{k-1}^{(i)}(z)}\,,
\]
\[
(\kappa^{(i)}_{1}\cdots\kappa^{(i)}_{k-1})^{2}=\frac{c_{1}^{(i)}\,c_{k-1}^{(i)}}{c_{k}^{(i)}}\,.
\]
Therefore, using  (\ref{defphi}), we get
\begin{equation}\label{simple}
\frac{\widetilde{F}_{k}^{(i)}(z)}{\widetilde{F}_{k-1}^{(i)}(z)(\kappa^{(i)}_{1}\cdots\kappa^{(i)}_{k-1})^{2}}=
\varphi_{k-1}^{(i)}(z)\,.
\end{equation}

Taking into consideration (\ref{ratioepsilons2})-(\ref{simple}),
we conclude that
\[
\lim_{\mathbf{n}\in\Lambda}\lambda_{\mathbf{n}}^*
\frac{\varepsilon_{\mathbf{n}_0,k-1}}
{\widetilde{\varepsilon}_{\mathbf{n},k-1}}
\frac{K_{\mathbf{n}_0,k-1}^2}{\widetilde{K}_{{\bf
n},k-1}^{2}}\frac{Q_{\mathbf{n}_0,k-1}(z)}
{\widetilde{Q}_{\mathbf{n},k-1}(z)}\frac{\widetilde{Q}_{\mathbf{n},k}(z)}
{Q_{\mathbf{n}_{0},k}(z)}(p_{k}\cdots p_{m})(z)=
\]
\[
c\,\Xi_{k}\prod_{i=1}^{k-1}\sign{p_{i},\supp{\sigma_{i}}}
\prod_{i=1}^{m-1}(\varphi_{k-1}^{(i)}(z))^{\deg(p_{i+1}\cdots
p_{m})} \lim_{\mathbf{n}\in\Lambda}\frac{Q_{\mathbf{n},k-1}(z)}
{\widetilde{Q}_{\mathbf{n},k-1}(z)}
\]
\begin{equation}\label{eqnuevo}
\times(p_{k}\cdots
p_{m})(z)\lim_{\mathbf{n}\in\Lambda}\frac{K_{\mathbf{n},k-1}^2}{\widetilde{K}_{{\bf
n},k-1}^{2}}\frac{\widetilde{Q}_{\mathbf{n},k}(z)}{Q_{\mathbf{n},k}(z)}\,,
\end{equation}
provided that the limits on the right hand side exist.

In Theorem, \ref{teo:1} we proved (\ref{asintrelativtilde}) for
$k=1$. Assume that $k=2$. Equations (\ref{eqasintotica}) and
(\ref{eqnuevo}) yield
\[
\lim_{\mathbf{n}\in\Lambda}\frac{K_{\mathbf{n},1}^2}{\widetilde{K}_{{\bf
n},1}^{2}}\frac{\widetilde{Q}_{\mathbf{n},2}(z)}{Q_{\mathbf{n},2}(z)}
\]
\[
=\frac{p_{\Lambda}(\varphi_{1}(z)){\mathcal{F}}(z;p_1,\ldots,p_m)
}{c\,\Xi_{2}\,\sign{p_{1},\supp{\sigma_{1}}}(p_{2}\cdots
p_{m})(z)\prod_{i=1}^{m-1}(\varphi_{1}^{(i)}(z))^{\deg(p_{i+1}\cdots
p_{m})}}\,,
\]
uniformly on compact subsets of $\overline{\mathbb{C}} \setminus
\supp{\sigma_{2}}$. Using (\ref{eqcasokigual2}), we have
\[
\frac{p_{\Lambda}(\varphi_{1}(z)) }{c\,\Xi_{2}\,(p_{2}\cdots
p_{m})(z)\,\prod_{i=1}^{m-1}(\varphi_{1}^{(i)}(z))^{\deg(p_{i+1}\cdots
p_{m})}}=G_{1}(z;p_{1},\ldots,p_{m})\,.
\]
Consequently,
\[
\lim_{\mathbf{n}\in\Lambda}\frac{K_{\mathbf{n},1}^2}{\widetilde{K}_{{\bf
n},1}^{2}}\frac{\widetilde{Q}_{\mathbf{n},2}(z)}{Q_{\mathbf{n},2}(z)}=
\frac{{\mathcal{F}}(z;p_1,\ldots,p_m)\,G_{1}(z;p_{1},\ldots,p_{m})}
{\sign{p_{1},\supp{\sigma_{1}}}} \,.
\]
Evaluating at infinity, we obtain
(${\mathcal{F}}(\infty;p_1,\ldots,p_m)=1$)
\[
\lim_{\mathbf{n}\in\Lambda}\frac{K_{\mathbf{n},1}^2}{\widetilde{K}_{{\bf
n},1}^{2}}=\frac{G_{1}(\infty;p_{1},\ldots,p_{m})}{\sign{p_{1},\supp{\sigma_{1}}}}\,.
\]
Therefore, (\ref{asintrelativKcuad}) and (\ref{relationFk}) are
satisfied for $k=1$, since $G_{0}={\mathcal{F}}$.

Define the functions
\[
\mathcal{F}_{k}(z;p_1,\ldots,p_m):=
\lim_{n\in\Lambda}\frac{\widetilde{Q}_{\mathbf{n},k}(z)}{Q_{\mathbf{n},k}(z)}
\]
provided the limit exists. From (\ref{eqcasokmayortres}) it
follows that for any $k\geq 3$,
\[
\frac{p_{\Lambda}(\delta\varphi_{k-1}(z))
}{c\,\Xi_{k}\,(p_{k}\cdots
p_{m})(z)\,\prod_{i=1}^{m-1}(\varphi_{k-1}^{(i)}(z))^{\deg(p_{i+1}\cdots
p_{m})}}=G_{k-1}(z;p_{1},\ldots,p_{m})\,.
\]
As a consequence, using (\ref{eqnuevo}), we obtain that for any
$k\geq 3$,
\[
\lim_{\mathbf{n}\in\Lambda}\frac{K_{\mathbf{n},k-1}^2}{\widetilde{K}_{{\bf
n},k-1}^{2}}\frac{\widetilde{Q}_{\mathbf{n},k}(z)}{Q_{\mathbf{n},k}(z)}=
\frac{\mathcal{F}_{k-1}(z;p_1,\ldots,p_m)\,G_{k-1}(z;p_{1},\ldots,p_{m})}{\prod_{i=1}^{k-1}
\sign{p_{i},\supp{\sigma_{i}}}}\,.
\]
Therefore, using an induction process one proves
(\ref{asintrelativtilde})-(\ref{relationFk}). \hfill $\Box$

\begin{co} \label{cor6.1}
Let $S = {\mathcal{N}}^*(\sigma_1,\ldots,\sigma_m)$. Consider the
perturbed Nikishin system
$\mathcal{N}(\frac{p_1}{q_1}\sigma_1,\ldots,\frac{p_m}{q_m}\sigma_m)$,
where $p_k,q_k$ denote relatively prime polynomials with real
coefficients whose zeros lie in $\mathbb{C} \setminus \cup_{k=1}^m
\Delta_k$. Let $\Lambda \subset {\mathbb{Z}}^{m}_+(\circledast)$
be a sequence of distinct multi-indices such that for all ${\bf n}
\in \Lambda$, $n_{1} - n_{m} \leq C,$ where $C$ is a constant. Let
$\widetilde{Q}_{{\bf n},k}, 1\leq k\leq m,$ be the monic
polynomials of degree $N_{{\bf n},k}$ whose simple zeros are
located at the points where $\widetilde{\Psi}_{n,k-1}$ vanishes on
$\Delta_{k}$, where $\widetilde{\Psi}_{{\bf n},k}, 0\leq k\leq m,$
denote the second type functions defined in
$(\ref{defnPsinktilde})$, with $p_{k}$ replaced by $p_{k}/q_{k}$.
Let $\widetilde{K}_{{\bf n},k}, 1\leq k\leq m$ be the constants
defined in $(\ref{defnKnktilde})$, with $p_{k}$ replaced by
$p_{k}/q_{k}$. Then, for each $k\in\{1,\ldots,m\}$,
\begin{equation}\label{asintrelativtilderat}
\lim_{\mathbf{n}\in\Lambda}\frac{\widetilde{Q}_{\mathbf{n},k}(z)}
{Q_{\mathbf{n},k}(z)}=\frac{\mathcal{F}_{k}(z; p_1 ,\ldots,p_m)}
{\mathcal{F}_{k}(z; q_1 ,\ldots,q_m)}\,, \quad K \subset
\overline{\mathbb{C}} \setminus \supp{\sigma_{k}}\,,
\end{equation}
and
\begin{equation}\label{asintrelativKcuadrat}
\lim_{\mathbf{n}\in\Lambda}\frac{\widetilde{K}_{\mathbf{n},k}^{2}}
{K_{{\bf n},k}^{2}}=\prod_{i=1}^{k}
\sign{p_{i}/q_{i},\supp{\sigma_{i}}}
\frac{G_{k}(\infty;q_{1},\ldots,q_{m})}{G_{k}(\infty;p_{1},\ldots,p_{m})}\,.
\end{equation}
\end{co}
{\bf Proof.} By $Q_{\mathbf{n},k}^{*}$ denote polynomials
associated with to the auxiliary Nikishin system ${\mathcal{N}}
(\sigma_1/{q_{1}},\ldots,
 {\sigma_m}/{q_{m}} )$, corresponding to the indices ${\bf n},k$. On account of Theorem
 \ref{theorelativepol}, we have that
\[
\lim_{\mathbf{n}\in\Lambda}\frac{\widetilde{Q}_{\mathbf{n},k}(z)}
{Q_{\mathbf{n},k}^{*}(z)}=\mathcal{F}_{k}(z; p_1 ,\ldots,p_m)\,,
\quad K \subset \overline{\mathbb{C}} \setminus
\supp{\sigma_{k}}\,.
\]
and
\[
\lim_{\mathbf{n}\in\Lambda}\frac{Q_{\mathbf{n},k}(z)}
{Q_{\mathbf{n},k}^{*}(z)}=\mathcal{F}_{k}(z; q_1 ,\ldots,q_m)\,,
\quad K \subset \overline{\mathbb{C}} \setminus
\supp{\sigma_{k}}\,.
\]
Therefore, (\ref{asintrelativtilderat}) is obtained. Using the
same idea, (\ref{asintrelativKcuadrat}) follows  from
(\ref{asintrelativKcuad}). \hfill $\Box$

\begin{rmk} Theorem \ref{relasympstf} and Corollary \ref{cor5.1}
allow to define polynomials $\widetilde{Q}_{{\bf n},k}, \\
k=1,\ldots,m,$ in the case when $p_k,q_k$ have complex
coefficients as those monic polynomials which carry the zeros of
$\widetilde{\Psi}_{{\bf n},k-1}$ which lie in $\mathbb{C}\setminus
\Delta_{k-1}$. For such polynomials $\widetilde{Q}_{{\bf n},k}$,
results analogous to those expressed in Theorem
\ref{theorelativepol} and Corollary \ref{cor6.1} can be proved.
\end{rmk}

{\bf Acknowledgments.} Both authors acknowledge support from
grants MTM 2006-13000-C03-02  of Min. de Ciencia y Tecnolog\'{\i}a
and CCG 06--UC3M/ESP--0690 of Universidad Carlos III de
Madrid-Comunidad de Madrid.


\begin{thebibliography}{99}

\bibitem{kn:Sasha}
A.I. Aptekarev. {\em Strong asymptotics of multiply orthogonal
polynomials for Nikishin systems}. Sbornik: Math. {\bf 190}
(1999), 631-669.

\bibitem{AptLopRoc} A.I. Aptekarev, G. L\' opez Lagomasino, and I.A.
Rocha. {\em Ratio Asymptotic of Hermite-Pade orthogonal
polynomials for Nikishin systems}. Mat. Sb.  {\bf 196} (2005),
1089-1107.

\bibitem{DolBer} D. Barrios Rolan\'ia, B. de la Calle Ysern, and G. L\'{o}pez
Lagomasino. {\em Ratio and relative asymptotic of polynomials
orthogonal with respect to varying Denisov-type measures}. J. of
Approx. Theory {\bf 139} (2006), 223-256.


\bibitem{kn:B-G}
B. de la Calle Ysern and G. L\'{o}pez Lagomasino. {\em Weak
Convergence of varying measures and Hermite-Pad\'{e} orthogonal
polynomials}. {\sl Constr. Approx.}  {\bf 15} (1999), 553-575.

\bibitem{Den} S.A. Denisov. {\em On Rakhmanov's theorem for Jacobi
matrices}. Proc. Amer. Math. Soc. {\bf 132} (2004), 847-852.

\bibitem{kn:Gonchar}
A.A. Gonchar, E.A. Rakhmanov, and V.N. Sorokin. {\em
Hermite-Pad\'{e} for systems of Markov-type functions}. Sbornik:
Math. {\bf 188} (1997), 671-696.

\bibitem{AG} A. L\'{o}pez Garc\'ia and G. L\'opez Lagomasino. {\em Ratio
asymptotic of Hermite-Pad\'{e} orthogonal polynomials for Nikishin
systems. II}. submitted.

\bibitem{kn:Nikishin}
E.M. Nikishin. {\em On simultaneous Pad\'{e} approximations}.
Math. USSR Sb. {\bf 41} (1982), 409-426.

\bibitem{kn:Rak3}
E.A. Rakhmanov. {\em On asymptotic properties of orthogonal
polynomials on the unit circle with weights not satisfying
Szeg\H{o}'s condition}.  Math. USSR Sb. {\bf 58} (1987), 149-167.

\end{thebibliography}
\end{document}